\documentclass{article}
\usepackage{amsthm}
\usepackage{amssymb}
\usepackage{amsmath}
\usepackage[mathscr]{eucal}
\usepackage{epsfig}
\usepackage{graphicx}
\usepackage{epic}
\usepackage{tikz}
\usepackage{xcolor}
\usepackage{bm}
\usepackage{caption}
\usepackage{subcaption}
\usepackage{multirow}
\usepackage[hyphens]{url}
\newtheorem{mytheorem}{\bf Theorem}[section]

\newtheorem{mylemma}[mytheorem]{\bf Lemma}
\newtheorem{corollary}[mytheorem]{\bf Corollary}
\newtheorem{remark}[mytheorem]{Remark}

\newcommand{\Tau}{\mathcal{T}}
\newcommand{\Mp}{\mathcal{P}}
\newcommand{\bx}{\bm{x}}
\newcommand{\bb}{\bm{b}}
\newcommand{\Eps}{\mathcal{E}}

\title{Robust  domain decomposition methods for high-contrast multiscale problems on irregular domains with virtual element discretizations}

\author{Juan G. Calvo\thanks{Centro de Investigación en Matemática Pura y Aplicada, Escuela de Matemática, Universidad de Costa Rica, Montes de Oca, San José, Costa Rica, 11501 \tt  juan.calvo@ucr.ac.cr.}
\and   Juan Galvis\thanks{Departamento de Matem\'aticas, Universidad Nacional de Colombia, Carrera 45 No. 26-85, Edificio Uriel Guti\'errez, Bogot\'a D.C., Colombia, \tt jcgalvisa@unal.edu.co.}
}

\date{May 2023}

\begin{document}

\maketitle

\begin{abstract}
Our research focuses on the development of domain decomposition preconditioners tailored for second-order elliptic partial differential equations. Our approach addresses two major challenges simultaneously: i) effectively handling coefficients with high-contrast and multiscale properties, and ii) accommodating irregular domains in the original problem, the coarse mesh, and the subdomain partition. The robustness of our preconditioners is crucial for real-world applications, such as the efficient and accurate modeling of subsurface flow in porous media and other important domains.

The core of our method lies in the construction of a suitable partition of unity functions and coarse spaces utilizing local spectral information. Leveraging these components, we implement a two-level additive Schwarz preconditioner. We demonstrate that the condition number of the preconditioned systems is bounded with a bound that is independent of the contrast. Our claims are further substantiated through selected numerical experiments, which confirm the robustness of our preconditioners.
\end{abstract}


\tableofcontents

\section{Introduction}
Nowadays, there is a wide variety of iterative methods for the solution of practical problems that involve the numerical approximation of solutions of partial differential equations (PDEs) in two or three dimensions. A complete review of all of these methodologies is impossible in a short manuscript, and we will focus on Domain Decomposition Methods (DDM), which have been effectively used for many problems in computational mathematics and engineering in the last century with tremendous success; see \cite{tw, dolean2015introduction} and references therein. These methods are based on ``divide and conquer" ideas and parallel computing paradigms in order to efficiently compute solutions of large ill-conditioned linear systems arising from discretizations of PDEs. Suppose that the solution of a PDE is computed by using a Finite Element Method (FEM) or a Virtual Element Method (VEM) with a fine enough mesh; see \cite{BrennerScott} and \cite{beirao}, respectively. Generally speaking, an iterative method based on DDM uses the efficient solution of small problems in order to obtain approximations of the solution of the original large problem. The ``small'' problems are either obtained by posing differential equations on subdomains or by using very coarse mesh resolutions. Then, in each iteration, solutions to these smaller problems are combined in order to converge to the solution of the original problem in a few iterations; see \cite{tw}.

Nevertheless, there are some particular cases where the classical recipes of DDM do not result in efficient computations, and some modifications have to be introduced in the methodology in order to recover efficient iterations. As it is usual, in order to propose the correct changes, a deep understanding of the analysis of the method is required which guides the modifications that shall be introduced. Usually, the efficiency is lost due to the need for many iterations in order to obtain approximations to some prescribed tolerance. Theoretically, this means that the operator driving the iteration is still very ill-conditioned. Another common reason is that the usual construction of partition of unity functions (needed in the analysis or in the computations) are inadequate due to irregular meshes or other artifacts. We remark two important examples that have been very active research areas in the last decades: 

\begin{itemize}
\item[(i)]  the case of PDEs with high-contrast and multiscale coefficients (\cite{ge09_2,Efen_GVass_11,egw10,galvis2017overlapping}), and 
\item[(ii)] the case of the presence of irregular domains in the original problem or in the subdomains (\cite{D&K&W_DDLessRegSubd,Widlund_2008,D&W_AltCoarseSpace,calvo2018approximation,calvo2019overlapping,calvo2020new}). 
\end{itemize}

In the case of (i), the presence of high-contrast multiscale coefficients in the PDE brings additional difficulties that are, in nature and at first glance, different from the challenges generated by (ii) that are related to the construction of a partition of unity with given overlapping supports with irregular boundaries. Addressing simultaneously both situations is crucial for real-world applications, such as modeling efficiently and accurately subsurface flow in porous media, among other important applications.

In order to fix ideas, we focus on the following elliptic equation with heterogeneous coefficients. Given a polygonal domain $D\subset \mathbb{R}^2$ and $f:D \to \mathbb{R}$, we consider the problem: Find
$u:D\to \mathbb{R}$ such that
\begin{equation}\label{eq:mainproblem}
-\mbox{div}(\kappa\nabla u)=f \mbox{ in } D,
\end{equation}
with Dirichlet boundary conditions $u=0$ on $\partial D$. 
The coefficient  $\kappa=\kappa(x)$
represents the permeability of the porous media $D$. 
In particular, we consider the case where $\kappa$ has significant local variations across $D$ at
different scales (multiscale structure). We assume that $\kappa$ is a high contrast coefficient, where the contrast is defined as
$$ \eta=\max_{x\in \overline{D}}\kappa(x)/\min_{x\in \overline{D}}\kappa(x).$$

The variational formulation for problem \eqref{eq:mainproblem} is: Find $u\in H_0^1(D)$ such that
\begin{equation}\label{eq:problem}
a(u,v)=f(v) \quad \mbox{ for all } v\in H_0^1(D),
\end{equation}
where the bilinear form $a(\cdot,\cdot)$ and the linear functional $f$ are defined by
\begin{equation*}
a(u,v)=\int_D 
\kappa(x)\nabla u(x)\cdot\nabla v(x) \ dx
\quad \mbox{ for all } u,v\in H_0^1(D),
\end{equation*}
and 
\[
f(v)=\int_Df(x)v(x)\ dx \quad \mbox{ for all } v\in H_0^1(D).
\]

In this work, we consider usual triangulations of the domain $D$, composed by triangles or squares (FEM), or general polygonal meshes such as hexagonal and Voronoi meshes (VEM). We refer to $\mathcal{T}^h$ as the fine-scale mesh, where $h$ is the size of a typical element, and we assume it is fine enough to resolve all the variations of the coefficient $\kappa$. We denote by $V^h= V^h(D)$ the finite or virtual element space chosen, and we write $V^h_0=V^h_0(D)=V^h(D)\cap H_0^1(D)$. The Galerkin formulation of (\ref{eq:problem}) is to
find $u\in V^h_0(D)$ such that 
\begin{equation}\label{eq:galerkin}
a(u,v)=f(v) \quad \mbox{ for all } v\in V^h_0(D),
\end{equation}
or in matrix form
\begin{equation}\label{eq:matrix}
Au=b,
\end{equation}
where, for all $u,v\in V^h(D)$, we have 
\[
u^TAv=\int_D\kappa \nabla u \cdot \nabla v 
\quad \mbox{ and }\quad   v^Tb=\int_D fv.
\]
For simplicity, we denote a finite or virtual element function and its vector of coordinates on $V_h$ with the same symbol.

It is well known that the performance of iterative methods for the solution of \eqref{eq:matrix} with a high-contrast multiscale coefficient depends on $\eta$ and on the local variations of $\kappa$ across $D$; see \cite{ge09_1,ge09_1reduceddim, ge09_2, abreu2019convergence} and references therein for the case of piecewise linear elements. This is also true for VEM, as observed in \cite{calvo2018approximation}. 
We remark that the mesh has to be fine enough to resolve the variations of the coefficient $\kappa$ in order to obtain good approximation results for classical FEM and VEM. Under these conditions, the discretization leads to the solution of very large (sparse) ill-conditioned problems (with the condition number scaling with $h^{-2}\eta$). Therefore, the performance of solvers  depends on $\eta$ and local variations of $\kappa$ across $D$. Before a solution to this issue was proposed in \cite{ge09_1,ge09_2}, some partial solutions were proposed in several works; see for instance, \cite{aarnes, Graham1, ge09_1}. 

In order to obtain an efficient good approximation of the solution $u$, two possible strategies are:

\begin{itemize}
\item  Choose $h$ sufficiently small and implement an iterative method with preconditioner $M^{-1}$ to solve $M^{-1}Au=M^{-1}b$, such that the condition number of $M^{-1}A$ is small and bounded independently of physical parameters such as $\eta$ and the multiscale structure of $\kappa$.
\item  Solve a smaller dimensional linear system: 
introduce a coarse-scale mesh $\mathcal{T}^H$ with {$H>h$} such that computations of solutions can be done efficiently\footnote{The coarse mesh does not necessarily resolve all the variations of $\kappa$.}. This usually involves the construction of  a downscaling operator (from the coarse-scale  to fine-scale $v_0\mapsto v$) and an upscaling operator $R_0$ (from fine-scale to coarse-scale, $v\mapsto v_0$) (or similar operators). Using these operators, the linear system $Au=b$ becomes  a coarse linear system $A_0 u_0=b_0$ so that  $R_0u_0$ or functionals of it can be computed. The main goal of this approach is to obtain a sub-grid capturing such that $|| u- R_0u_0||$  is small; see for instance \cite{eh09,egh12}.
\end{itemize}
In this paper, we consider the first of these two approaches and we focus on the construction of a DD preconditioner $M^{-1}$. We construct a two-level additive Schwarz DDM. The first level consists of adding up solutions to local problems related to an overlapping decomposition $\{ D_i'\}_{i=1}^{N_S}$ of the domain $D$ where the problem is posed. The second level needs the construction of a coarse approximation related to a coarse mesh $\mathcal{T}^H$. In particular, for high-contrast  multiscale problems, the construction of the coarse space needs to be adapted to the variations of the coefficient either by adapting the coarse mesh (and/or the non-overlapping decomposition) to the discontinuities of the coefficient or by adapting the construction of the coarse basis functions to the variations of the coefficient. For high-contrast multiscale coefficients, the former is not always practical; see
\cite{ge09_1reduceddim,ge09_2} and related works. 

Moreover and also very important in real-world applications, in many cases it is impossible to construct an initial domain decomposition into non-overlapping subdomains in such a way that we obtain regular interfaces between them. This is the case when the fine mesh is not regular or structured in the original problem and the application does not suggest a natural partition. It is usual to use partitioning software such as ParMETIS
\footnote{ParMETIS (Parallel Graph Partitioning and Fill-reducing Matrix Ordering) is an MPI-based parallel library that implements a variety of algorithms for partitioning unstructured graphs, meshes, and for computing fill-reducing orderings of sparse matrices. More information can be found in \url{https://www3.cs.stonybrook.edu/~algorith/implement/ParMetis/implement.shtml}.} 
that, loosely speaking, groups fine-scale elements, and very few requirements can be imposed to obtain smooth interfaces. There are also important applications where the permeability coefficient has very complicated multiscale variations. In these cases, the partition algorithm, even if it does not take into account the variation of the coefficients, it generates partitions with irregular interfaces. For such irregular decompositions, classical constructions of coarse basis functions and a partition of unity do not satisfy the small gradient restriction needed in the analysis and/or computations. Therefore, the coarse basis functions and/or partition of unity functions need to be adapted to the irregular shape of the subdomains; see \cite{Widlund_2008, calvo2018virtual, calvo2019overlapping} and related works. 

From the previous discussion, high-contrast coefficients and irregular subdomains bring some important challenges in order to design efficient iterative methods based on DDM. To the best of the author's knowledge, these two cases are usually considered separately in all the revised literature; i.e., either there is an assumption on smoothness or boundedness of the coefficient, or on the smoothness of the boundaries of the subdomains (Liptchitz assumption). In this paper, we do not assume either of these two cases. Our purpose is to design a DD iterative method whose performance deteriorates minimally with respect to the simultaneous presence of (i) multiscale variations in the coefficient; (ii) large jumps and high-contrast; and (iii) irregular subdomains and interfaces.

In order to obtain a robust DDM, the construction of appropriate partition of unity functions and coarse spaces are the key ingredients to focus on. Therefore, we center our attention on the construction of partition of unity functions and coarse spaces that can handle both difficulties at the same time.

The rest of the paper is organized as follows. In Section \ref{sec:vem} we describe virtual spaces, which are required for the construction of the partition of unity and the discretization of problem \eqref{eq:galerkin} in the presence of general polygonal elements. In Section \ref{sec:OSDD} we review two-level overlapping Schwarz methods.  We then study how to obtain a stable decomposition in Section \ref{sec:stableDD}; we review the construction of coarse spaces for high-contrast multiscale problems and for irregular subdomains, and combine these two constructions in order to obtain a robust method that works for our problem. Finally, in Section \ref{sec:exp} we present some numerical results and in Section \ref{sec:conc} we include some conclusions and final remarks.

\section{Virtual element spaces} \label{sec:vem}
VEM for problems posed in $H^1(D)$ were introduced in \cite{MR2997471}, while some practical aspects were presented in \cite{beirao}. This method can handle general polygonal elements, which is a natural choice to consider in the presence of irregular subdomains. We briefly discuss the method and refer to \cite{MR2997471} for a more detailed explanation of the method. For numerical implementation details we refer to \cite{Sutton,MR4544725} for the cases $k=1$ and $k\geq 2$, respectively, where $k$ denotes the degree of polynomial spaces that are considered. We remark that we can avoid discrete harmonic extensions by appropriately projecting virtual functions to polynomial spaces of degree $k\geq 2$, which allows us to approximate our partition of unity functions for irregular subdomains; see \cite{calvo2018approximation}.

Given an integer $k\geq 1$, the local virtual space of degree $k$ is defined for every element $E\in\Tau_h$ as
\begin{equation*} 
V_k(E) := \lbrace v\in H^1(E): v|_{\partial E}\in C^0(\partial E),\ v|_{e}\in\Mp_k(e)\ \forall {\rm edge }\ e,\ \Delta v\in \Mp_{k-2}(E)\rbrace,
\end{equation*}
where $\Mp_{-1}(E) = \lbrace 0\rbrace$. Given $v \in V_k(E)$, its local degrees of freedom (dof) can be chosen as:
\begin{enumerate}
\item The value of $v$ at each vertex of $E$.
\item The value of $v$ at the $k-1$ internal Gauss-Lobatto quadrature nodes of every edge of $E$.
\item The moments of $v$ up to order $k-2$ given by $$\frac{1}{\vert E\vert} \int_{E} v p\quad \forall p\in \Mp_{k-2}(E);$$ 
\end{enumerate}
see \cite[Equation 2.4]{MR2997471}. The total number of dof is then $N_{\rm dof}^E := kN_V^E+k(k-1)/2$, where $N_V^E$ is the number of vertices of $E$. We remark that there are equivalent sets of dof that can be used. Define the functional 
$${\rm dof}_i^E : V_k(E)\rightarrow \mathbb{R},$$ 
where ${\rm dof}_i ^E(v)$ is the $i$-th degree of freedom of $v$. The canonical local basis is then defined as the set $\lbrace \varphi_i^E\rbrace_{i=1}^{N_{\rm dof}^E}$ such that ${\rm dof}_i^E(\varphi_j^E) = \delta_{ij}$.  Define the projector $P_0: V_k(E)\rightarrow \Mp_0(E)$ by 
\begin{align*}
P_0 v &= \dfrac{1}{N^E_V} \sum_{i=1}^{N^E_V} v(\bx_i),\quad \text{for }\ k=1;\\
P_0 v &=\dfrac{1}{\vert E\vert} \displaystyle \int_E v,\quad \text{for }\ k\geq 2.
\end{align*}

The $L^2$ projection of the gradient $\Pi_{E,k}^\nabla: V_k(E)\rightarrow \Mp_k(E)$ is defined as follows. Let $w=\Pi_{E,k}^\nabla v$ be the polynomial of degree at most $k$ such that $$
a^E(w,p) = a^E(v,p)\quad \forall\ p\in \Mp_k(E)\quad \text{and}\quad P_0 v = P_0 w,$$
where $\displaystyle a^E(u,v) = \int_E \nabla u \cdot \nabla v$. The local bilinear form is defined as
\begin{equation*}
a_h^E(u,v) = \int_E \nabla\Pi_{E,k}^\nabla u \cdot \nabla \Pi_{E,k}^\nabla v + \sum_{r=1}^{N^E_{\text{dof}}} \text{dof}_r((I-\Pi_{E,k}^\nabla)u) \text{dof}_r((I-\Pi_{E,k}^\nabla)v);
\end{equation*}
see \cite[Section 3.3]{beirao} for further details. Therefore, the entries corresponding to the local stiffness matrix are given by
\begin{equation*}
(A^E_h)_{i,j} = a_h^E(\varphi_j,\varphi_i).
\end{equation*}

Globally, the virtual element space $V_h\subset H_0^1(D)$ is defined as $$V_h(D) = \lbrace v\in H_0^1(D): v|_{E} \in V_k(E)\ \text{ for all\ } E\rbrace.$$ 
The global bilinear form is then obtained by assembling the local bilinear forms:
\begin{equation*}
a_h(u,v) = \sum_{E }a_h^E(u,v).
\end{equation*}
 
For the local right hand side, different approaches have been proposed; see, e.g., \cite[Section 6]{beirao} and references therein. We consider the $L^2(E)-$projection operator $\Pi^0_k: V_k(E)\rightarrow \Mp_k(E)$. Then, the local contribution for the load term can be defined by
\begin{equation*} 
(\bb_E^h)_i = \int_E f\ \Pi^0_k \varphi_i \quad (i=1,\ldots,N^E_\text{dof}).
\end{equation*}
Standard optimal estimates can be obtained for the solution of the linear system; see, e.g., \cite{MR3073346}. 

\section{A brief review of robust overlapping methods} \label{sec:OSDD}
In this section, we present some important aspects related to the design and analysis of two-level overlapping DDM. This material is clearly focused on what we need to present for our work and does not pretend to be an exhaustive literature review.

\subsection{The classical additive two-level method} 

Given a non-overlapping decomposition $\{ D_i\}_{i=1}^{N_S}$ of the domain $D$ (that is conforming with respect to the fine-mesh $\mathcal{T}^h$), we construct an overlapping decomposition $\{ D_i'\}_{i=1}^{N_S}$ and denote by 
$\delta$ the maximum diameter of the overlap among subdomains. For instance, we can add layers of fine-mesh elements around each $D_i$. For this case, we denote by $\delta_i$ the maximum width of the region $D_i'\setminus D_i$ and define $\delta = \max_{1\leq i\leq N_S} \delta_i$. Let $V^i_0(D_i')$ be the set of finite or virtual element functions with support in $D_i'$. We also denote by $R_i^T:V^i_0(D_i')\to V^h$ the zero extension operator.

Let $A_i$ be the  Dirichlet matrix corresponding to the overlapping subdomain $D_i'$ defined by 
\begin{equation*}
v^TA_iw=a(R_i^Tv,R_i^Tw) \quad \mbox{ for all } v,w\in V^i_0(D'_i),
\end{equation*}
and $i\in\{1,\dots,N\}$; see \cite{tw} and references therein.  The classical one-level additive  method solves
\begin{equation*}
 M_1^{-1}Au=M_1^{-1}b,
\end{equation*}
 with
 \begin{equation*}
 M^{-1}_1=\sum_{i=1}^{N_S} R_i^T{A}_i^{-1} R_i. 
\end{equation*}
The application of the one-level preconditioner  involves
solving  $N_S$ local problems in each iteration. We have the bound $$\mbox{Cond}(M^{-1}_1A)\preceq \left(1+{1}/{\delta H}\right);$$ for high-contrast multiscale problems it is known that $\displaystyle C \asymp \eta$; see \cite{ge09_1,ge09_2}.

We will introduce a set of coarse basis functions $\{\Phi_i\}_{i=1}^{N_c}$ defined on the coarse triangulation $\mathcal{T}^H$, where $N_c$ is the number of coarse basis functions. Define then the coarse space 
$$V_0=\mbox{span}\{\Phi_i\}_{i=1}^{N_c}.$$  The coarse-scale matrix  is $$A_0=R_0AR_0^T,$$ where 
$$R_0^T=[\Phi_1,\dots,\Phi_{N_c}].$$
The two-level preconditioner
uses the coarse space and it is defined by 
\begin{equation}\label{eq:Minvsub2}
M^{-1}_2=  R_0^TA_0^{-1}R_0+
\sum_{j=1}^{N_S} R_j^T ({A}_j)^{-1} R_j= R_0^TA_0^{-1}R_0+M_1^{-1}. \end{equation}
It is known that 
\begin{equation*}
\mbox{Cond}(M^{-1}_2A)\preceq
\eta
\left(1+{H}/{\delta}\right).
\end{equation*}

The classical two-level method is robust with respect to the number of subdomains but it is not robust with respect to $\eta$.
The condition number  estimates use a Poincar\'e
inequality and a small overlap trick; see, for example \cite{tw}. Without  the small overlap trick, it holds that
$\mbox{Cond}(M^{-1}A)\preceq \eta(1+H^2/\delta^2)$ as is recalled in \cite{ge09_1,ge09_1reduceddim}. To our knowledge, there is not a high-contrast independent small overlap trick. 

\subsection{High-contrast multiscale coefficients case}
\begin{figure}
\centering
\psfig{figure=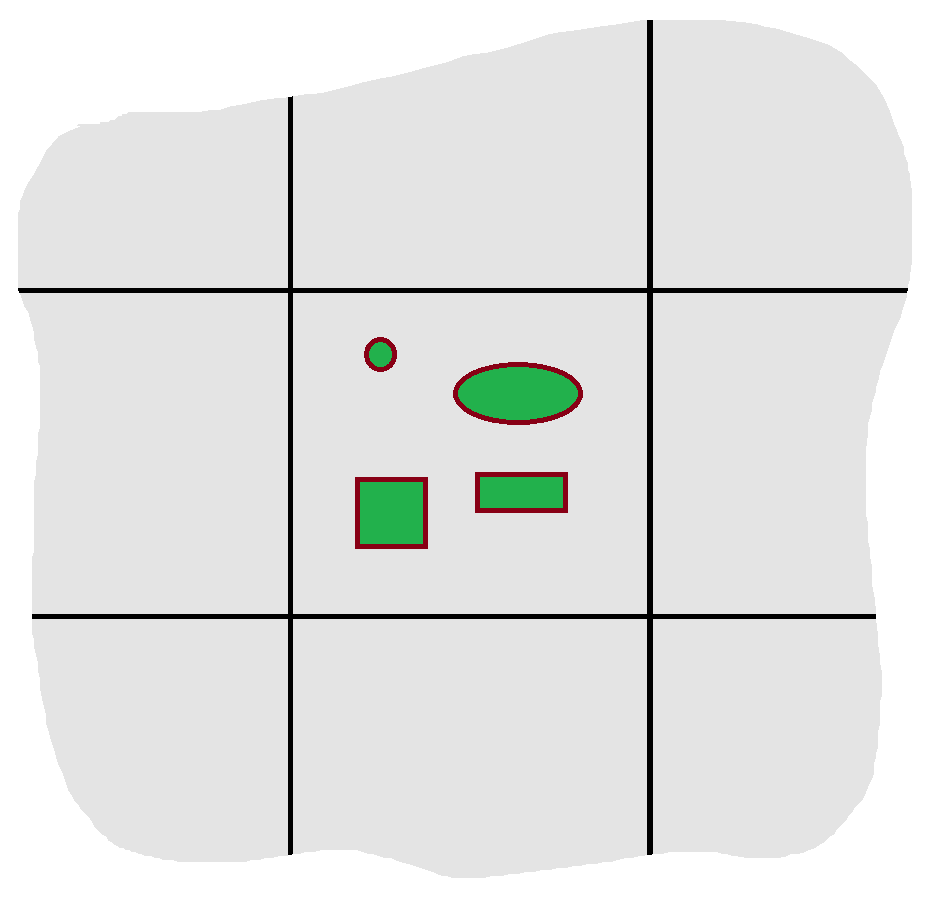,height=3cm,width=3cm,angle=0}
\psfig{figure=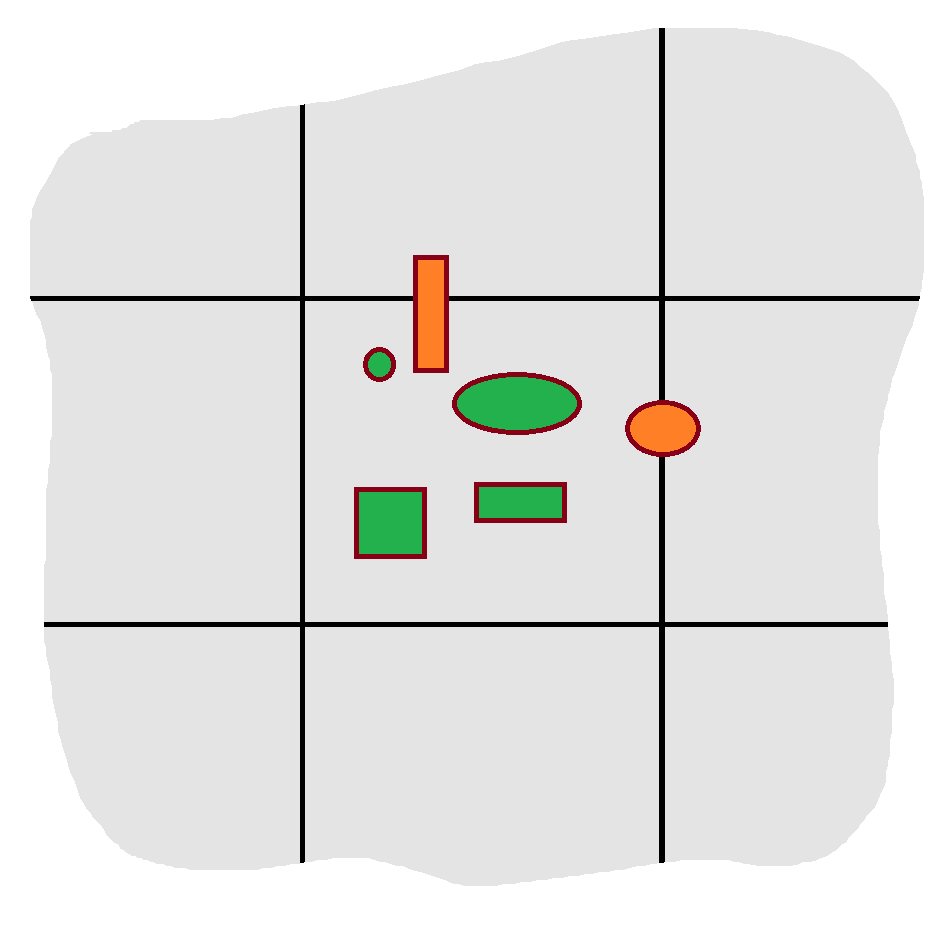,height=3cm,width=3cm,angle=0}
\psfig{figure=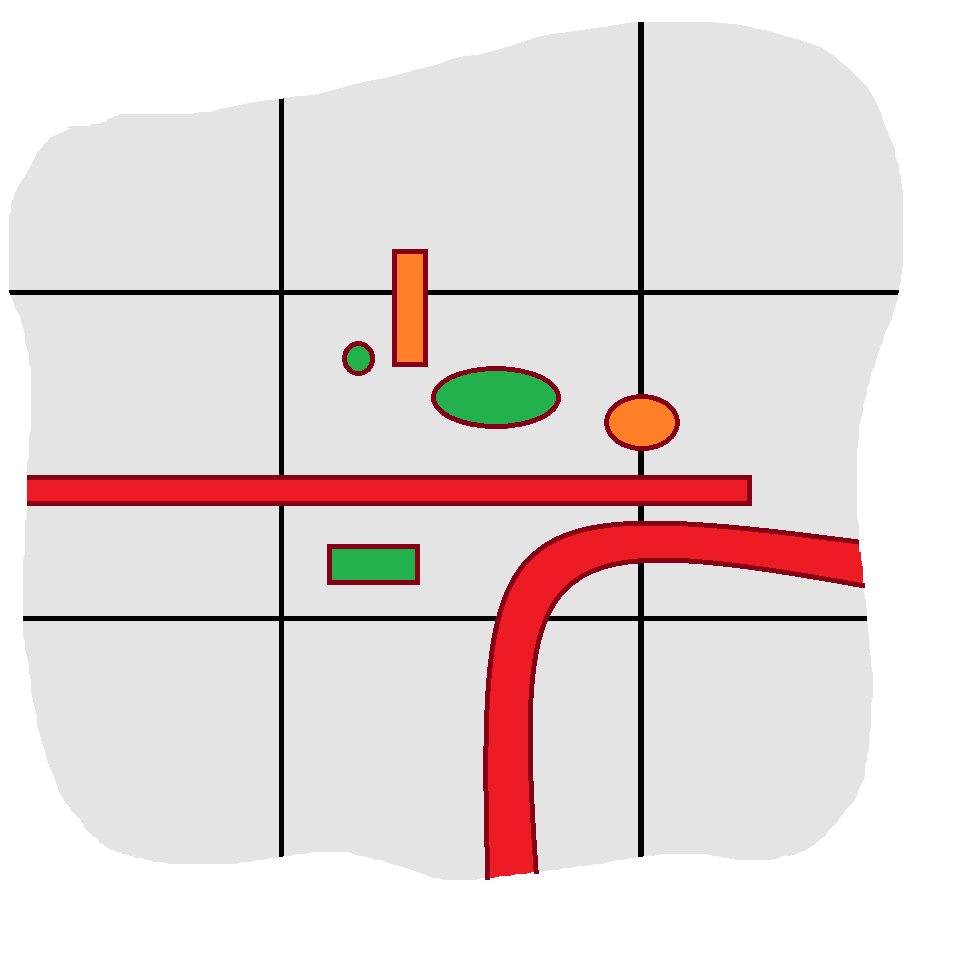,height=3cm,width=3cm,angle=0}
\caption{Examples of multiscale  subdomains  with interior high-contrast inclusions in green (left), boundary inclusions in orange (center), and long channels in red (right). We have $\kappa=1$ in the gray background, and $\kappa = \eta \gg 1$ inside the channels and inclusions.}\label{inclusions}
\end{figure}

For the diffusion coefficient $\kappa$ in problem 
\eqref{eq:mainproblem}, we assume that there exists $\kappa_{\min}$ and $\kappa_{\max}$ such that 
\[
0< \kappa_{\min}\leq \kappa(x)\leq \kappa_{\max}\]
for all  $x\in \overline{D}$. For $\Omega \subset D$, consider the contrast of $\kappa$ restricted to $\Omega$ defined by
\[
\eta_\Omega =\frac{\max_{x\in \Omega} \kappa(x)}{\min_{x\in \Omega}\kappa(x)}.
\]
We assume that $\kappa\geq 1$ is a piecewise constant coefficient with background equal to 1; this is, we assume that there is a family of open disjoint subsets $\{Q_i\}_{i=0}^{N_R}$ with
$\overline{D}=\bigcup_{i=0}^{N_R} \overline{Q}_i$, such that 
\[
\kappa(x)=\begin{cases}
\kappa_i \geq 1 & x\in \overline{Q}_i, \\
1, & x\in Q_0,
\end{cases}
\]
and $\partial D \subset \partial Q_0$. The piecewise-smooth coefficient case is similar. We call $Q_0$ the background region and the sets $Q_i$ are referred to as high-value inclusions or high-value channels; see Figure \ref{inclusions} for an illustration. 

Denote by $\mathcal{S}^H(D)$ the coarse-scale triangulation and/or partition of the domain into non-overlapping subdomains, with $H\gg h$. We then say that the coefficient $\kappa$ has a multiscale high-contrast structure with respect to the triangulation $\mathcal{T}^h$ and the partition $\mathcal{S}^H$ if 
\begin{itemize}
    \item[(a)] for every $e\in \mathcal{T}^h$ we have that $\kappa|_e$ is constant; and
    \item[(b)] there exists $\alpha\in (0,1]$ and  $\mathcal{S}_0^H\subset \mathcal{S}^H$ such that  for all $E\in \mathcal{S}^H_0$ we have that 
    $\eta\geq \eta_E\geq \alpha\eta$.
\end{itemize}
We refer to subdomain $E$ as a high-contrast subdomain. Note that the label ``multiscale coefficient'' refers to the fact that the coefficient varies at multiple scales and therefore there is a large quantity of high-contrast subdomains scattered everywhere around the whole domain.

There are several works addressing the performance of classical DDM for high-contrast problems. Many of these works
considered simplified multiscale structures\footnote{These works usually assume some alignment between the coefficient heterogeneities and the initial non-overlapping decomposition.}; see, e.g., \cite{tw} for some works by O. Widlund and his collaborators. We also mention the work by Sarkis and his collaborators, where they introduced the assumption of quasi-monotonicity \cite{MR1367653}. Sarkis also introduced the idea of using ``extra'' or additional basis functions, as well as techniques that construct the coarse spaces using the overlapping decomposition (and not related to a coarse mesh); see \cite{MR2099424}. 

Hou and Aarnes \cite{aarnes}, and Scheichl and  Graham \cite{Graham1}, started a systematic study of the performance of classical overlapping DDM for high-contrast problems. In their works, they used coarse spaces constructed using a coarse grid and special basis functions from the family of  multiscale finite element
methods. These  authors designed two-level DDM that were robust with respect to $\eta$, for special multiscale structures aligned with the coarse mesh. For example, the mentioned cases presented good performance for the first two cases shown in Figure \ref{inclusions}, but not for the third one that includes long high-contrast channels (long high-conductivity regions that touch the boundaries of more than two subdomains). 

None of the results available in the literature (before the method introduced in \cite{ge09_1,ge09_1reduceddim}) were robust for a general coefficient not-aligned with the construction of the coarse space (not aligned either with the non-overlapping decomposition or the coarse mesh if any); i.e., the condition number of the resulting preconditioner depended strongly on $\eta$ for general multiscale coefficients. Many developments and extensions of the methodology in \cite{ge09_1,ge09_1reduceddim, abreu2019convergence} have been proposed by many authors, including extensions not mentioned here such as extensions or adaptations to other DD preconditioners such as the non-overlapping methods, among others. We only mention the work in \cite{eglw11}, where the method was generalized to the case of any positive bilinear form that can be assembled from local bilinear forms. 

\subsection{Irregular subdomains and partition of unity}\label{sec:ispu}
A bounded domain $\Omega\subset\mathbb{R}^2$ is a uniform domain (in the sense of Peter Jones \cite{Jones1981}) if there exists a constant $C_U(\Omega)>0$ such that for any pair of points ${a}$, ${b}$ in the closure of $\Omega$, there is a curve $\gamma(t):[0,l]\rightarrow\Omega$, parametrized by arc length, with $\gamma(0)={a}$, $\gamma(l)={b}$ and with
\begin{equation}
l\leq C_U\vert{a}-{b}\vert,\notag
\end{equation}
\begin{equation}
\min(\vert\gamma(t)-{a}\vert,\vert\gamma (t)-{b}\vert)\leq C_U \text{dist}(\gamma(t),\partial\Omega).\notag
\end{equation}
For a rectangular domain we have $C_U\geq L_1/L_2$ where $L_1$, $L_2$ are the height and width of the domain. Thus, the constant $C_U$ can be large if the subdomain has a large aspect ratio.


Jones domains form the largest class of finitely connected domains for which an extension theorem holds in two dimensions; see \cite[Theorem 4]{Jones1981}. DDM for irregular subdomains and scalar elliptic problems in the plane for Jones domains started in \cite{D&K&W_DDLessRegSubd, D&W_AltCoarseSpace, Widlund_2008}; there are similar studies for problems posed in $H(\mbox{rot})$ and $H(\mbox{div})$. Previous studies constructed coarse basis functions by defining the degrees of freedom on the boundary of the subdomains and then extending these values to the interior of each subdomain (nodes in the fine mesh) with discrete harmonic extensions. Bounds are similar to previous studies, but the construction of $R_0^T$ involves the solution of a linear system for each subdomain with a size equal to the number of interior nodes (fine mesh) for each subdomain. A variation was considered in \cite{calvo2018approximation,calvo2018virtual}, where appropriate polynomial projections are considered for the interior of the subdomains, reducing considerably the size of the problems that have to be solved. We remark that these coarse functions add up to 1, so they are a natural candidate for partitions of unity for irregular subdomains.

We assume that the polygonal subdomains $\{D_j\}_{j=1}^{N}$ are uniform domains and we denote the lowest virtual space on the coarse mesh by 
$$V_H=V_H(D) = \lbrace v\in H_0^1(D): v|_{D_i} \in V_1(D_i)\ \text{ for all\ } D_i\rbrace.$$ 
For a subdomain vertex $y_i$, consider the set $\omega_i$ given by the union of subdomains that share $y_i$. We will define a coarse (virtual) function $\chi_i\in V_H$, with compact support on $\omega_i$, such that $\{\chi_i\}$ will be a partition of unity for the overlapping partition $\{\omega_i\}$. For that, we follow similar ideas as in \cite{calvo2018virtual}, where we take the columns of $R_0^T$ as functions for the partition of unity. 

\begin{figure}
     \centering
\includegraphics[width=\textwidth]{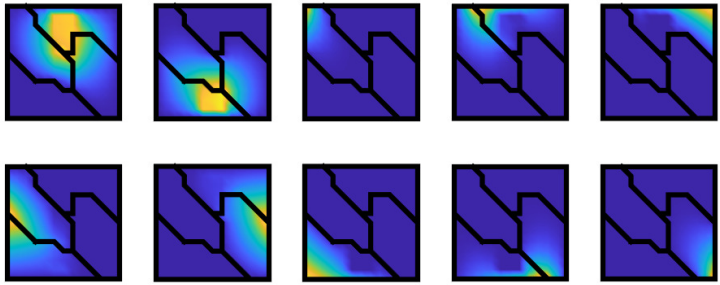}
     \caption{Partition of unity $\{\chi_i\}$ for 4 METIS subdomains obtained from a triangular mesh, based on discrete harmonic extensions. We have two functions related to internal subdomain vertices, and eight functions related to vertices on the boundary of the domain. Nodal values vary from 0 (blue) to 1 (yellow).}
        \label{fig:partition}
\end{figure}

For each subdomain vertex $y_i$, we define $\chi_i \in V_H$ by choosing appropriately its degrees of freedom; see Figure \ref{fig:partition} for the case of METIS subdomains. First, we set $\chi_i(y) = 0$ for all the subdomain vertices $y$, except at $y_i$ where $\chi_i(y_i) = 1$. Second, we set the degrees of freedom related to the nodal values on each subdomain edge $\Eps$ (defined as $\Eps=\Eps^{ij}=\overline{D}_i\cap \overline{D}_j$). If $y_i$ is not an endpoint of $\mathcal{E}$, then $\chi_i$ vanishes on that edge. If $\mathcal{E}$ has endpoints $y_i$ and $y_j$, let ${d}_\mathcal{E}$ be the unit vector with direction from $y_i$ to $y_j$. Consider any node $\widetilde{y} \in \mathcal{E}$. If $0\leq (\widetilde{y} - y_{j})\cdot {d}_{\Eps} \leq \vert y_i -  y_j\vert$, we then set 
\begin{equation*}
\chi_i(\widetilde{{y}}) = \dfrac{(\widetilde{y} - y_{j})\cdot {d}_{\Eps}}{\vert y_i - y_j\vert}.
\end{equation*}

It is clear that $\chi_i(y_i) = 1$, $\chi_i(y_j) = 0\ (j\neq i)$, and that the function varies linearly in the direction of ${d}_\mathcal{E}$ for such nodes. If $(\widetilde{y} - y_{j})\cdot {d}_{\Eps} < 0$ or $(\widetilde{y} - y_{j})\cdot {d}_{\Eps} > \vert y_i -  y_j\vert$, we then set $\chi_i(\widetilde{y}) = 0$  or $\chi_i(\widetilde{y}) = 1$, respectively. In this way, we define all the degrees of freedom of $\chi_i \in V_H$. By construction, it is clear that $0\leq \chi_i\leq 1$, $\sum_{i} \chi_i=1$, ${\rm supp}(\chi_i)\subset {\omega_i}$ and, if we assume typical hypothesis on the coarse mesh for VEM, it holds that $\|\nabla \chi_i\|_{\infty} \leq C/H$, for  some constant $C$ that depends on the logarithm of the maximum ratio of distances between consecutive subdomain nodes. 

Instead of using discrete harmonic extensions to define values in the interior nodes of subdomains, we can project virtual functions to polynomial spaces by just knowing the values on the interface; see the construction of the operator $R_0^T$ in \cite{calvo2018virtual}. This approach saves computational time, since the minimization problem requires to solve just a $\mathcal{O}(k^2)$ linear system, where $k$ is the degree of the polynomial space. Even though these functions no longer satisfy to be bounded by 0 and 1, they add up to one and satisfy similar bounds as the discrete harmonic extensions, and numerically, they provide competitive results; see \cite{calvo2018approximation}.

\section{Stable decomposition under high-contrast and irregular subdomains} \label{sec:stableDD}

A main tool in obtaining condition number bounds is the construction of a stable decomposition of a global field; see \cite{tw} and \cite{galvis2021condition} for an abstract interpretation. Hereafter we denote $v_j\in V_j=V_h(D_j')$ and its extension by zero to $D$ by the same symbol. 

A decomposition for a global field $v\in V=V_h(D)$ is written as 
\begin{equation*}
 v=v_0+\sum_{j=1}^{N_S} v_j,
\end{equation*}
where $v_0\in V_0$ corresponds to the coarse part of the decomposition and $v_j\in V_j$ ($j=1,\dots, N$) correspond to the local parts of the decomposition.

The decomposition is stable in the sense that there exists $C_0>0$ such that
\begin{equation}\label{eq:estable_decomposition}
\int_{D} \kappa |\nabla v_0|^2+\sum_{j=1}^{N_S}
\int_{D_j'}\kappa |\nabla v_j|^2 \leq C_0^2\int_{D}\kappa |\nabla v|^2.
\end{equation}
Equation \eqref{eq:estable_decomposition} implies the bound (see \cite{tw})
$$\mbox{cond}(M_2^{-1}A)\leq c(\mathcal{T}^h,\mathcal{T}^H) C_0^2,$$
where $c(\mathcal{T}^h,\mathcal{T}^H)$ is a constant that depends on the fine and coarse grid configurations; e.g., the numbers of coarse blocks that intersect a subdomain. Moreover, the existence of a suitable coarse interpolation $I_0:V\to V_0=\mbox{span}\{\Phi\}$ allows the use of the coarse space $V_0$ in the stable decomposition above; see \cite{aarnes,tw,Graham_Scheichl_07}. Usually such stable decomposition is constructed as explained next.

\subsection{Coarse part of the stable decomposition}
We denote by $\{y_i\}_{i=1}^{N_v}$ the subdomain vertices of the coarse mesh
$\mathcal{T}^H$ and define the neighborhood of the 
node $y_i$ by
\begin{equation*}
\omega_i=\bigcup\{ K_j\in\mathcal{T}^H; ~~~ y_i\in \overline{K}_j\}
\end{equation*}
and the neighborhood of the coarse element $K$ by
\begin{equation}\label{eq:def:omegaK}
\omega_K=\bigcup \{ \omega_j\in\mathcal{T}^H; ~~~ y_j\in 
\overline{K}_j\}.
\end{equation}

For the coarse part of the stable decomposition, we
introduce a partition of unity  $\{\chi_i\}$ subordinated to the coarse mesh (supp $\chi_i\subset \omega_i$). As remarked in Section \ref{sec:ispu} and \cite{calvo2018approximation,calvo2018virtual,Jones1981}, special attention and methods have to be call out when constructing such a partition of unity in the case of coarse blocks with irregular boundaries. 

We  begin by restricting the global field $v$ to  $\omega_i$.
For each coarse node neighborhood $\omega_i$, we
identify a local field that will contribute to the
coarse space, denoted by
{$  I_0^{\omega_i}v $},  in such a way that the
coarse space will be defined as\footnote{We mention that this identification is similar to the identification of multi-continua averages in some recent versions of  the generalized multiscale finite element method; see \cite{chung2017coupling} and related works.}  $$V_0= \mbox{Span}\{ \chi_i  I_0^{\omega_i}v\}.$$
In classical two-level domain decomposition methods methods  $I_0^{\omega_i}v$ is the average of $v$ in
$\omega_i$, but in \cite{ge09_1,ge09_1reduceddim} more general examples for $I_0^{\omega_i}$ where introduced related to projection in local eigenspaces of a carefully chosen eigenvalue problem.

We assemble a global coarse field as
\begin{equation*}
v_0=I_0v=\sum_{i=1}^{N_S} \chi_i (I_0^{\omega_i}v ).
\end{equation*} 
The operator $I_0$ shall be constructed to be stable in the $H^1$ norm and have some approximation properties in the $L^2$ norm. Note that in each block we have 
\begin{equation*}
v-v_0=\sum_{x_i\in K} \chi_i (v-  I_0^{\omega_i}v ).    
\end{equation*}

\subsection{Local part of the stable decomposition}
For the local parts of the stable decomposition, we introduce a partition of unity  $\{\xi_j\}$ subordinated to the overlapping decomposition (supp $\xi_j\subset  D_j'$); see Section \ref{sec:ispu} and \cite{calvo2018approximation,calvo2018virtual} for the construction of such a partition of unity when the subdomains have irregular boundaries. 

The local part of the stable decomposition is defined by
\[
v_j=I^h(\xi_j (v-v_0)),
\]
where $I^h$ is the fine-scale interpolation. To bound the energy of  $v_j$ we observe that, 
\begin{eqnarray}\int_{D_j' } \kappa|\nabla v_j|^2  &\preceq &
\int_{D_j^\prime}\kappa
|\nabla I^h(\xi_j(v-v_0))|^2 
\label{eq:energyofvj} 
\preceq \int_{D_j^\prime}\kappa
|\nabla (\xi_j(v-v_0))|^2, \nonumber 
\end{eqnarray}
and  in each coarse-block $K$ we have, 
\begin{eqnarray*}
\int_{K}\kappa
|\nabla (\xi_j(v-v_0))|^2 
 &\preceq& \int_{K} \kappa\left|\nabla \left( \xi_j\left(\sum_{x_i\in K} \chi_i (v-  I_0^{\omega_i}v )\right)\right)\right|^2 \label{eq:gradpu1}\\
  &\preceq&  \sum_{x_i\in K} \int_{K} \kappa\left|\nabla \left( \xi_j\chi_i (v-  I_0^{\omega_i}v )\right)\right|^2. 
   \label{eq:gradpu2}
\end{eqnarray*}

Adding up over $K$ and using 
$$\widehat{D}_j
=\bigcup \{ \overline{K} \quad  | \quad K \cap D_j'\not= \emptyset \},
$$
we obtain,
\begin{eqnarray}
\int_{D_j' } \kappa|\nabla v_j|^2 &\preceq& 
\sum_{K \cap D_j'\not= \emptyset } \int_{K} \kappa\left|\nabla \left( \xi_j\chi_i (v-  I_0^{\omega_i}v )\right)\right|^2\nonumber \\
&\preceq & \sum_{x_i\in  \widehat{D}_j}\int_{\omega_i }\kappa\left|\nabla \left( \xi_j\chi_i (v-  I_0^{\omega_i}v )\right)\right|^2. \label{eq:smaloverlapstop}
\end{eqnarray}
Note that the hidden constants above depend on the 
relation between the partition into overlapping subdomains and the coarse mesh and also on the connectivity of the coarse mesh. These constants are not essential for the main point of our paper which is the irregularity of the subdomains and coarse blocks, and the presence of multiscale discontinuous high-contrast coefficients. To obtain a stable decomposition, we would like to bound the last term by $C \int_{\Omega_j } \kappa |\nabla v|^2$, where $D_j'\subset \Omega_j $.

\begin{remark}\normalfont
The roles of the partition of unity related to the coarse mesh and the partition of unity related to the subdomain decomposition are essential in the analysis. However, explicit gradient bounds are not needed for the implementation. If we use the analysis presented in Section \ref{secabs}, we will not need explicit bounds of gradients of the partitions of unity functions. In this paper, for the sake of simplicity and concreteness, we present the construction of appropriate multiscale coarse spaces making use of explicit bounds of the gradients of the partition of unity functions. The point of avoiding explicit bounds of gradients of partition of unity becomes important when analyzing and implementing multilevel methods with three of more levels. See \cite{efendiev2012multiscalepg,egvdd20} for a multilevel extension of the methods presented here where the construction of additional levels is done in an algebraic way.
\end{remark}

At this point we have identified the following cases that are usually considered in the reviewed literature: 
\begin{itemize}
    \item {\bf Coarse mesh-based partition:} the subdomain  overlapping decomposition coincides with coarse node neighborhoods; i.e., $D_j'=\omega_j$.
 In this case, we can write a bound with $\xi$ instead of $\chi$ and
 replace $\nabla (\chi^2)$ by $\nabla \chi$ and therefore  we need to bound the term
 $$\sum_{x_i\in \overline{\omega}_j} \int_{\omega_j} \kappa|\nabla   \chi_i |^2| v-
I_0^{\omega_i}v  |^2.$$  In this case, the coarse-scale operator $I_0$, in addition to having the  $H^1$ stability and $L^2$ approximation properties, it also has approximation properties in the $H^1$ norm; see for instance 
\cite{egw10,abreu2019convergence} and related works. 
\item {\bf Overlap-based coarse space:} the coarse basis functions support coincides with the overlapping decomposition. Here also $\omega_i=D_j'$. In this case, there is no coarse mesh per se, and $V_0= \mbox{Span}\{ \xi_i  I_0^{\omega_i}v\}$. Therefore, the variation of coarse-scale basis functions occurs only on the union of overlaps; see \cite{spillane2014abstract} and related works.
\item {\bf General case:} The overlapping partition and the coarse scale partition are constructed independently.  This is the case considered in this article, emphasizing 
that both subdivisions, the coarse mesh and the subdomain decomposition, may have irregular boundaries (a common practical situation in real-world applications); see \cite{calvo2018approximation, calvo2019overlapping, calvo2020new}. For the general construction for the case of smooth subdomains and high-contrast bilinear forms, see \cite{eglw11}.
\end{itemize}

As noted in \cite{eglw11}, we can stop in inequality  \eqref{eq:smaloverlapstop} and proceed to obtain abstract bounds. Here, in order to fix ideas, we go a step further and bound, 
\begin{eqnarray*}
\int_{\omega_i }\kappa\left|\nabla \left( \xi_j\chi_i (v-  I_0^{\omega_i}v )\right)\right|^2 &\preceq &
 \int_{\omega_i } \kappa(\xi_j\chi_i)^2|\nabla (v-I_0^{\omega_i}v)|^2 \\&&
+{    \int_{\omega_i} \kappa|\nabla (\xi_j\chi_i)|^2| v-
I_0^{\omega_i}v  |^2},
\end{eqnarray*}
which allows us to concentrate on bounding the second term above. A literature review reveals that the main techniques to obtain such bounds are (\cite{galvis2017overlapping}): 
\begin{itemize}
\item {\bf Poincar\'e inequality:} Classical analysis for the case of bounded coefficient and smooth interfaces  uses a Poincar\'e inequality (sometimes) combined with the small overlap trick to obtain the required bound above; see \cite{Widlund_2008}. This analysis was extended for the case of irregular boundaries and bounded coefficients in \cite{D&K&W_DDLessRegSubd}. For high-contrast multiscale coefficients,
the resulting bound depends on the contrast $\eta$, in general. For the quasi-monotone coefficient, it can be obtained a bound that is independent of the contrast \cite{MR1367653} as well as for the case of {locally connected high-contrast regions} \cite{ge09_1}. In the general case is not possible to use this technique and obtain bounds independent of the contrast. 
\item 
{\bf $L^\infty$ estimates:} the idea is to use an $L^\infty$ estimate of the form
\begin{eqnarray*}
  \int_{\omega_i} {   \kappa|\nabla (\xi_j\chi_i)|^2}| v-
I_0^{\omega_i}v  |^2 
\preceq
   {  || \kappa|\nabla (\xi_j\chi_i)|^2 ||_\infty} {   \int_{\omega_i} | v- I_0^{\omega_i}v  |^2}.
\end{eqnarray*}
The idea in \cite{Graham1,aarnes}
was to construct  partitions  of unity such that the term $|| \kappa|\nabla (\xi_j\chi_i)|^2 ||_\infty$  is bounded independently of the contrast, and then to use classical Poincar\'e inequality estimates. 
Instead of minimizing the $L^\infty$, the authors in \cite{Graham1} intuitively  tried to minimize $\int_K \kappa|\nabla \chi_i|^2 $. In the light of 
\cite{D&K&W_DDLessRegSubd}, the same idea can be used for the case of irregular interfaces. 
This strategy works well when the multiscale structure of the coefficient is confined within the coarse blocks and the subdomains are regular. In the general case, the obtained bounds depend on the contrast.
\item {\bf Local generalized eigenvalue problems:} this is the technique used in this paper and will be explained next. This technique was introduced in \cite{ge09_1,ge09_1reduceddim} for high-contrast problems and smooth subdomains, and it is now a popular technique for obtaining an analysis of DDM that are required to be robust with respect to some parameters. Thanks to the works 
\cite{calvo2018approximation,calvo2018virtual,D&K&W_DDLessRegSubd} and related papers addressing the case or irregular boundary decompositions, we could successfully use the same method for the case of irregular subdomains. It is also worth to mention that \cite{ge09_1,ge09_1reduceddim}  motivated the introduction of a whole family of very flexible and robust multiscale methods now known as the Generalized Multiscale Finite Element Method; see \cite{egh12, galvis2020numerical,abreu2019convergence} and related works. 
\end{itemize}

\subsection{Local generalized eigenvalue problem} 
\begin{figure}[t]
     \centering
     \begin{subfigure}[b]{0.3\textwidth}
         \centering
         \includegraphics[width=\textwidth]{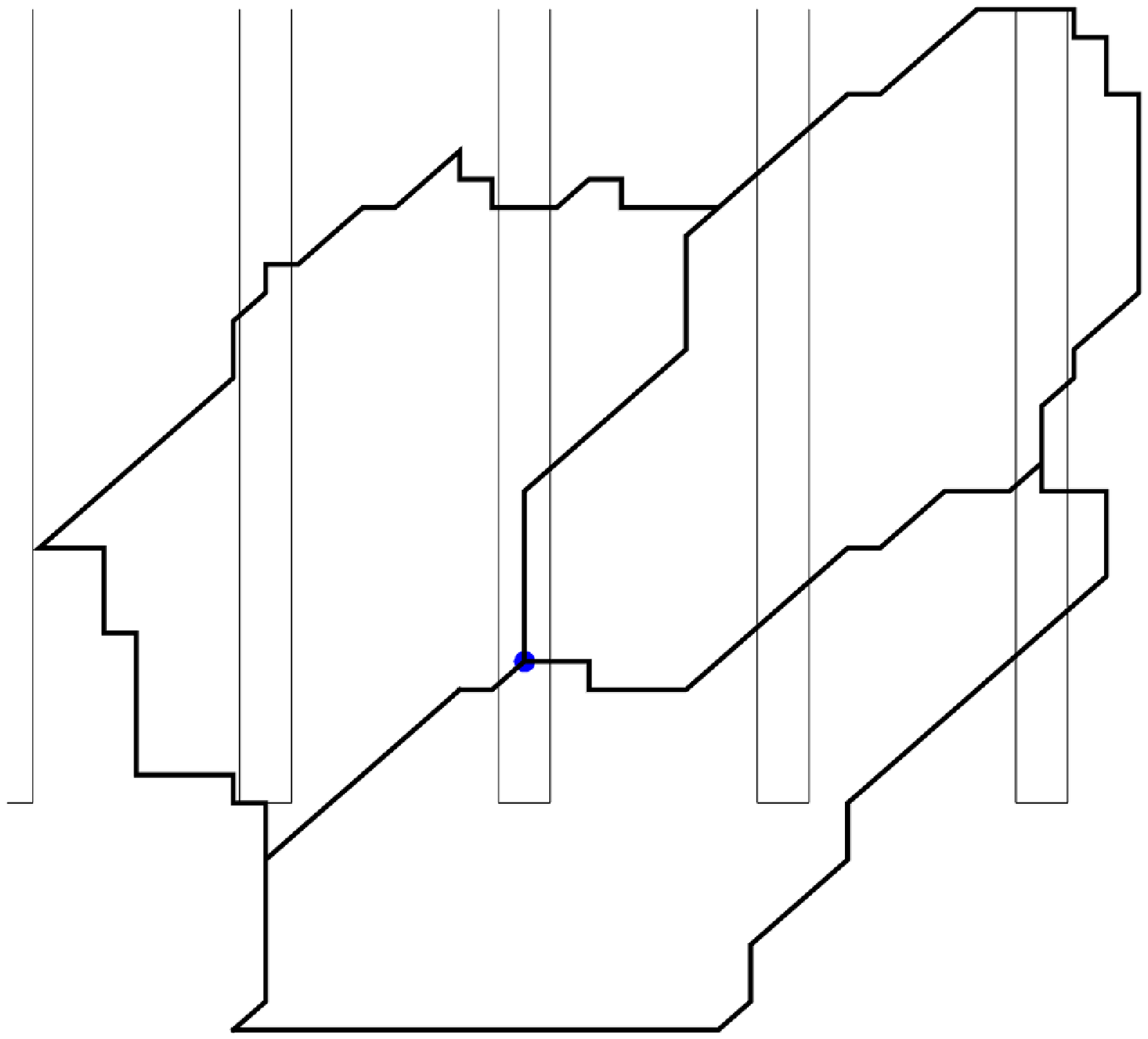}
     \end{subfigure}
     \hfill
     \begin{subfigure}[b]{0.6\textwidth}
         \centering
         \includegraphics[width=\textwidth]{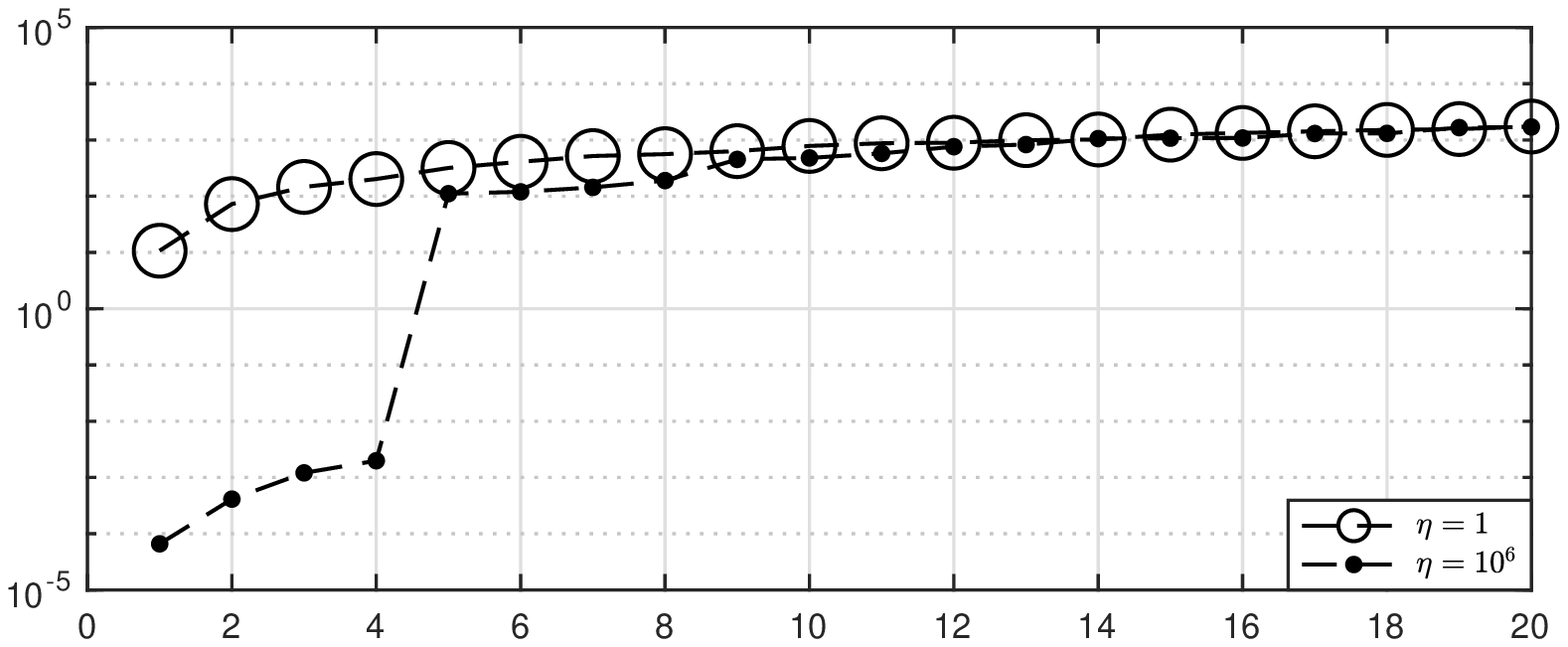}
     \end{subfigure}
    \caption{(left) A subdomain vertex with three METIS subdomains. The coefficient $\kappa$ is $\kappa = \eta$ inside the small rectangular channels, and $\kappa = 1$ in the background. 
    (right) Eigenvalue distribution for $\eta=1$  (circles) and $\eta=10^6$ (black dots). The effect of having a high contrast coefficient implies the addition of four eigenfunctions, associated to the four eigenvalues smaller than 1.}
        \label{fig:eigenvalues_eta0_6}
\end{figure}

We follow the proposals in \cite{ge09_1,ge09_2,eglw11} and related references. We can write
$${  \int_{\omega_i}  \kappa|\nabla (\xi_j\chi_i)|^2  |v-
I_0^{\omega_i}v |^2}
\preceq
{ \frac{1}{\delta^2 H^2} \int_{\omega_i}
\kappa|(v- I_0^{\omega_i}v )|^2}\preceq C
{  \int_{\omega_i} \kappa |\nabla v|^2},$$
where we need to justify the last inequality with constant independence of the contrast. The idea is then to consider the Rayleigh quotient,
$$\displaystyle \mathcal{Q}(v):=\frac{\int_{\omega_i} \kappa |\nabla v|^2}{\int_{\omega_i} \kappa| v  |^2}$$ with $v\in V_h(\omega_i).$
This quotient is related to an eigenvalue problem and we can define $ I_0^{\omega_i}v $ to be  the projection on low modes of this quotient  on $\omega_i$. The associated eigenproblem is given by 
\begin{equation}\label{eq:eigproblem}
-\mbox{div}(k(x)\nabla \psi_\ell^{\omega_i})=\lambda_\ell k(x)\psi_i^{\omega_i} \mbox{ in } \omega_i
\end{equation}
with homogeneous Neumann boundary conditions for floating subdomains and a mixed homogeneous Neumann-Dirichlet condition for subdomains that touch the boundary. It turns out that the low part of the spectrum can be written as

\begin{equation}\label{eq:ordering}
\lambda_1\leq \lambda_2\leq ...\leq \lambda_L<\lambda_{L+1}\leq 
\end{equation}
where   $\lambda_1 , ...,  \lambda_L$ are small, asymptotically vanishing eigenvalues and $\lambda_L$ can be bounded below independently of the contrast; see \cite{ge09_2,eglw11} for precise statements and proofs and Figure \ref{fig:eigenvalues_eta0_6} for a typical curve of the smallest eigenvalues. After identifying the local field $I_0^{\omega_i}v$,
we then define the coarse space as
\begin{equation*}
    V_0=\mbox{Span}\{ I^h\chi_i \psi_j^{\omega_i} \}= \mbox{Span}\{ \Phi_i \}.
\end{equation*}

We assume that the elements  
of $\mathcal{T}^h$ contained in $\omega_i$ form 
a triangulation of $\omega_i$. Define the matrix
$A^{\omega_i}$ corresponding to homogeneous Neumann problem by 
\begin{equation*}
v^TA^{\omega_i} w=\int_\Omega \kappa \nabla v \nabla w \quad \mbox{ for all } v,w\in {V}_h({\omega_i}),
\end{equation*}
and the {\it modified mass matrix} of same dimension $M^{\omega_i}$ by
\begin{equation}\label{eq:def:M^Ni}
v^T{M}^{\omega_i} w=\int_{\Omega}\widetilde{\kappa} v w \quad \mbox{ for all } v,w\in \widetilde{V}^h(\Omega),
\end{equation}
where $\widetilde{V}_h=V_h(\omega_i)$ if 
$\overline{\Omega}\cap \partial D=\emptyset$ and 
$\widetilde{V}_h=\{ v\in V_h(\omega_i) : v=0 \mbox{ on } 
\partial\omega_i\cap\partial D\}$ otherwise. 
Here $\widetilde{\kappa}$ in 
(\ref{eq:def:M^Ni}) is a weight 
derived from the high-contrast coefficient 
 $\kappa$ and  contains the relevant 
information we need for the construction of the coarse 
basis functions. In the above case $\widetilde{\kappa}=\kappa$ but other choices can also be used \cite{ge09_1reduceddim}.
We have then the generalized eigenvalue problem
\begin{equation}\label{eq:eigenvalueproblem}
A^{\omega_i}\psi={\lambda} {M}^{\omega_i}\psi.
\end{equation}

Let $n_h(\omega_i)$ denote the number of degrees of freedom 
in $\overline{\omega_i}$. Given any $v\in V^h(\omega_i)$ we can write 
\[
v=\sum_{\ell =1}^{n_h(\omega_i)}  \left(
v^TM^{\omega_i}\psi_\ell^{\omega_i} 
\right)\psi_\ell^{\omega_i} =
\sum_{\ell =1}^{n_h({\omega_i})}  \left(
\int_{{\omega_i}}\kappa v\psi_\ell^{\omega_i} 
\right)\psi_\ell^{\omega_i} 
\]
and compute 
\begin{equation}\label{eq:avv}
\int_{{\omega_i}} \kappa |\nabla v|^2=v^TA^{\omega_i} v=\sum_{\ell =1}^{n_h({\omega_i})}
 \left(
\int_{{\omega_i}}\kappa v \psi_\ell^{\omega_i}
\right)^2 \lambda_{\ell }^{\omega_i}
\end{equation}
and 
\begin{equation}\label{eq:mvv}
\int_{{\omega_i}}\kappa v^2=
v^T M^{\omega_i} v=
\sum_{\ell=1}^{n_h({\omega_i})}
 \left(
\int_{{\omega_i}}\kappa v\psi_\ell^{\omega_i} 
\right)^2 .
\end{equation}
Given an integer $L$ and $v\in V^h({\omega_i})$, we define 
\begin{equation}\label{eq:def:I-Omega-L}
I^{{\omega_i}}_{L}v=\sum_{\ell=1}^L \left(
\int_{{\omega_i}}\kappa v\psi_\ell^{\omega_i}
\right)\psi_\ell^{\omega_i}. 
\end{equation}
From (\ref{eq:ordering}),
(\ref{eq:avv}),  and (\ref{eq:mvv}) it is easy to prove the following inequality:
\begin{equation}\label{eq:truncation}
\int_{{\omega_i}} \kappa(v-I_L^{\omega_i} v)^2\leq 
\frac{1}{\lambda_{L+1}^{\omega_i}} a(v-I_L^{\omega_i} v,v-I_L^{\omega_i} v)
\leq 
\frac{1}{\lambda_{L+1}^\Omega} a(v,v).
\end{equation}

\begin{remark}\normalfont We make the following remark regarding \eqref{eq:eigenvalueproblem} and inequality \eqref{eq:truncation}: 
\begin{enumerate}
    \item When $L=1$, $\kappa=1$ (or 
    $\kappa$ is smooth and bounded)  and  $\partial \omega_i$ is  smooth (Lipschitz), we obtain the classical
Poincar\'e inequality since it can be verified that 
$\lambda_{2}^{\omega_i}\preceq \mbox{\normalfont diam}({\omega_i})^{-2}$, where $\mbox{diam}({\omega_i})$
is the diameter of ${\omega_i}$.
\item For the case of a subdomain 
$\Omega$ with smooth boundaries and high-contrast coefficient it can be verified that  if $L$ is large enough then $\lambda_{L+1}^{\Omega}$
is contrast independent. We refer to 
(\ref{eq:truncation}) in this case  as the contrast independent weighted Poincar\'e 
inequality. In fact, the number $L$ corresponds to the number of high-contrast channels and inclusions or the number of high-contrast channels depending on the initial partition of unity used. For details see \cite{Efendiev_Galvis_10} and related works.

\item Finally, for the case of irregular boundaries and high-contrast coefficients, we obtain a behavior similar to the previous case; see Figure         \ref{fig:eigenvalues_eta0_6} for an illustration. In this case also, the number of small (asymptotically vanishing) eigenvalues are related to the number of high-contrast regions. As before, for the case $\kappa=1$ we recover the Poincar\'e inequality for subdomains with irregular boundaries as in \cite{D&K&W_DDLessRegSubd}.
\end{enumerate}
\end{remark}

We note that $\{\omega_{i}\}_{y_i\in \mathcal{T}^H}$ 
is a covering of $\Omega$.
Let $\{\chi_i\}_{i=1}^{N_c}$ be a partition 
of unity subordinated to the covering
$\{\omega_{i}\}$ such that $\chi_i\in V^h(D)$ and
$|\nabla \chi_i|\leq \frac{C}{H}$, $i=1,\dots,N_c$, where the constant depends on the logarithm of the ratio of distances between consecutive subdomain vertices; see \cite[Lemma 2.7]{D&W_AltCoarseSpace} and Section\
\ref{sec:ispu}. Define the set of coarse 
basis functions
\begin{equation*}
\Phi_{i,\ell}=I^h(\chi_i\psi_\ell^{\omega_i}) \quad \mbox{ for } 1\leq i\leq N_c
\mbox{ and } 1\leq \ell \leq L_i,
\end{equation*}
where $I^h$ is the fine-scale nodal value 
interpolation and  $L_i$ is an integer number for each $i=1,\dots,N_c$.
Note that in this case, there are several basis functions per coarse node.
The number of basis functions per node is defined via 
the eigenvalue problem (\ref{eq:eigproblem}).
Denote by $V_0$ the \emph{local spectral multiscale} space   
\begin{equation*}
V_0=\mbox{span}\{ \Phi_{i,\ell}: 1\leq i\leq N_c
\mbox{ and } 1\leq \ell \leq L_i \}.
\end{equation*}
Define also the coarse interpolation $I_0:V^h(D)\to V_0$ by
\begin{equation*}
I_0v=\sum_{i=1}^{N_c} \sum_{\ell=1}^{L_i}
\left(\int_{\omega_i}\kappa v
\psi_\ell^{\omega_i} \right)
I^h(\chi_i\psi_{\ell}^{\omega_i})=
\sum_{i=1}^{N_c} I^h\left( (I^{\omega_i}_{L_i}v)\chi_i\right),
\end{equation*}
where $I^h$ is the fine-scale nodal value interpolation and $I^{\omega_i}_{L_i}$ is defined in (\ref{eq:def:I-Omega-L}). Since $\psi_\ell^{\omega_i}\in \widetilde{V}_h(\omega_i)$, we see that $I_0v$ satisfies the zero Dirichlet boundary condition on $\partial D$. Note that we have
\[
v-I_0v=\sum_{i=1}^{N_c} 
I^h\left(\chi_i(v-I^{\omega_i}_{L_i}v)\right).
\]

We have the following result. The proof (for the case of smooth subdomains) is presented in  
\cite{Efendiev_Galvis_10} and we present it for the sake of completeness. 
\begin{mylemma}\label{lem:coarse-projection}
For all coarse element $K$ we have 
\begin{equation}\label{eq:I0approx}
\int_{K}\kappa  (v-I_0v)^2\preceq 
\frac{1}{\lambda_{K,L+1}}
\int_{\omega_K} \kappa |\nabla v|^2,\end{equation}
\begin{equation}\label{eq:I0stab}
\int_{K}\kappa  |\nabla I_0v|^2\preceq 
\max\Big\{1,\frac{1}{H^2\lambda_{K,L+1}}\Big\}
\int_{\omega_K} \kappa |\nabla v|^2,
\end{equation}
where $\lambda_{K,L+1}=\min_{y_i\in K}\lambda_{L_i+1}^{\omega_i}$
and $\omega_K$ is the union of the elements that share common edge
with $K$ defined in (\ref{eq:def:omegaK}).
\end{mylemma}

\begin{proof} 
First we prove (\ref{eq:I0approx}). Using the fact that 
$\chi_i\leq 1$ we have 
\begin{eqnarray*}
\int_{K}\kappa  (v-I_0v)^2&\preceq& \sum_{y_i\in K}
\int_{K} \kappa
I^h(\chi_i(v-I^{\omega_i}_{L_i}v)^2 \\
&\preceq& 
 \sum_{y_i\in K}
\int_{K} \kappa (\chi_i(v-I^{\omega_i}_{L_i}v))^2 \\\
&\preceq&  \sum_{y_i\in K}
\int_{\omega_i} \kappa (v-I^{\omega_i}_{L_i}v)^2
\end{eqnarray*}
and using (\ref{eq:truncation}) with $\Omega=\omega_i$
to estimate the last term above, 
we obtain
\begin{eqnarray*}
\int_{K}\kappa  (v-I_0v)^2 &\preceq& \sum_{y_i\in K} 
\frac{1}{\lambda_{L+1}^{\omega_i}} \int_{\omega_i} \kappa
|\nabla v|^2\\
&\preceq& \max_{y_i\in K} 
\frac{1}{\lambda_{L+1}^{\omega_i}} \int_{\omega_K} \kappa
|\nabla v|^2.
\end{eqnarray*}
To prove the stability (\ref{eq:I0stab}) we note that $\sum_{y_i\in K} \nabla \chi_i=0$ in $K$,
and  then we can 
fix $y_j\in K$ and
write $\nabla\chi_j=
-\sum_{y_i\in K\setminus\{y_j\}} \nabla \chi_i$. 
We obtain,
\begin{eqnarray*}
\nabla  \sum_{y_i\in K}  (I^{\omega_i}_{L_i}v)\chi_i
&=&\sum_{y_i\in K}  
\nabla \chi_i  (I^{\omega_i}_{L_i}v) +
\sum_{y_i\in K} \chi_i  \nabla (I^{\omega_i}_{L_i}v) \\
&=&\sum_{y_i\in K\setminus\{y_j\}}  
(I^{\omega_i}_{L_i}v-I^{\omega_j}_{L_j}v)\nabla\chi_i  +
\sum_{y_i\in K} \chi_i\nabla(I^{\omega_i}_{L_i}v) 
\end{eqnarray*}
which gives the following bound valid on $K$, 
\begin{eqnarray} \label{eq:gradsumI0local}
|\nabla  \sum_{y_i\in K}  (I^{\omega_i}_{L_i}v)\chi_i|^2
&\preceq &\frac{1}{H^2}\sum_{y_i\in K\setminus\{y_j\}}  
(I^{\omega_i}_{L_i}v-I^{\omega_j}_{L_j}v)^2+
\sum_{y_i\in K} |\nabla(I^{\omega_i}_{L_i}v)|^2.
\end{eqnarray}
Since  $\sum_{y_i\in K}  (I^{\omega_i}_{L_i}v)\chi_i \in 
\mathbb{P}^3(K)$  we can use the stability of the fine grid interpolation and 
(\ref{eq:gradsumI0local}) to get
\begin{eqnarray}
\int_{K}\kappa |\nabla I_0v|^2 
&=& \int_{K} \kappa 
|\nabla  I^h(\sum_{y_i\in K}  (I^{\omega_i}_{L_i}v)\chi_i)|^2\nonumber\\
&\preceq& \int_{K} \kappa 
|\nabla  \sum_{y_i\in K}  (I^{\omega_i}_{L_i}v)\chi_i|^2
\nonumber \\ 
&\preceq& 
\sum_{y_i\in K} 
\frac{1}{H^2}\int_{K} \kappa (I^{\omega_i}_{L_i}v-I^{\omega_j}_{L_j}v)^2 
+\sum_{y_i\in K}\int_{K} \kappa |\nabla(I^{\omega_i}_{L_i}v)|^2.
\label{eq:twotermsgradI0}\end{eqnarray}
To bound the first term above we use 
(\ref{eq:truncation}) with $\Omega=\omega_i$ as follows, 
\begin{eqnarray}
\int_{K} \kappa (I^{\omega_i}_{L_i}v-I^{\omega_j}_{L_j}v)^2 
&\preceq & 
\int_{\omega_i} \kappa (v-I^{\omega_i}_{L_i}v)^2+
\int_{\omega_j} \kappa (v-I^{\omega_i}_{L_i}v)^2\nonumber\\
&\preceq & 
\frac{1}{\lambda_{L+1}^{\omega_i}} \int_{\omega_i} \kappa
|\nabla v|^2+
\frac{1}{\lambda_{L+1}^{\omega_j}} \int_{\omega_j} \kappa
|\nabla v|^2\nonumber\\
&\preceq &
\frac{1}{\lambda_{K,L+1}} \int_{\omega_K} \kappa
|\nabla v|^2. \label{eq:firstgradI0}
\end{eqnarray}
The second term in (\ref{eq:twotermsgradI0}) 
is estimated using (\ref{eq:avv}) and the orthogonality of the 
eigenvectors in the $A^{\omega_i}$ inner product
\begin{eqnarray}
\int_{K} \kappa |\nabla(I^{\omega_i}_{L_i}v)|^2\leq 
\int_{\omega_i} \kappa |\nabla(I^{\omega_i}_{L_i}v)|^2
\leq \int_{\omega_i} \kappa |\nabla v|^2\leq  
\int_{\omega_K} \kappa |\nabla v|^2.\label{eq:secondgradI0}
\end{eqnarray}
By combining  (\ref{eq:firstgradI0}), (\ref{eq:secondgradI0})
and (\ref{eq:twotermsgradI0}) 
we obtain (\ref{eq:I0stab}).\end{proof}

\begin{corollary}
Under the assumptions of Lemma \ref{lem:coarse-projection}, 
the condition number of the preconditioned operator 
$M_2^{-1}A$ with $M_2^{-1}$ defined in (\ref{eq:Minvsub2}) is 
of order
\[
\mbox{cond}(M_2^{-1}A)\preceq 
C_0^2\preceq\max\left\{1+ 
\frac{1}{\delta^2 \lambda_{L+1}},1+
\frac{1}{H^2\lambda_{L+1}}\right\}, 
\]
where 
$\displaystyle \lambda_{L+1}=\min_{1\leq i \leq N_c} \lambda_{L_i+1}^{\omega_i}$.

\end{corollary}

\begin{remark}
Assume that $\kappa(x)=1$ for all $x\in D$ and 
$L_i=1$ for all $i=1,\dots,N_c$.  
Then (\ref{eq:I0approx}) and (\ref{eq:I0stab}) become 
\[
\int_{K}\kappa  (v-I_0v)^2\preceq 
H^2
\int_{\omega_K} \kappa |\nabla v|^2,\ \ \ 
\int_{K}\kappa  |\nabla I_0v|^2\preceq 
\int_{\omega_K} \kappa |\nabla v|^2,
\]
respectively, 
because  $\lambda_{K,L+1}=
\max_{y_i\in K}\lambda_{L_i+1}^{\omega_i}\asymp H^{-2}$.
\end{remark}

\begin{corollary}Assume the coefficient $\kappa$ is a two-valued coefficient with value $1$ in the background and $\eta$ inside inclusions and high-contrast channels. 
If we select $L_i$, the number of high-contrast dependent eigenvalues (that correspond to the number of connected components of the region with high-contrast value), the condition number of the preconditioned operator $M_2^{-1}A$ with $M_2^{-1}$ defined in (\ref{eq:Minvsub2}) 
is of order
$
\mbox{cond}(M_2^{-1}A)\preceq C(1+ \frac{H^2}{\delta^2}),
$ where $C$ is independent of the contrast and the mesh size.
\end{corollary}

\subsection{Eigenvalue problem with a multiscale partition of unity.}
Instead of the argument presented earlier, we can include the gradient of the partition of unity in the bounds (somehow similar to the ideas of $L^\infty$ bounds). We then need the following 
\begin{equation*}
{   \int_{\omega_i}|v-
I_0^{\omega_i}v |^2}
=
{  \frac{1}{H^2} \int_{\omega_i}
\widetilde{\kappa}|v- I_0^{\omega_i}v )|^2}
\preceq
{   \int_{\omega_i} \kappa |\nabla v|^2},
\end{equation*}
where we have introduced
\begin{equation*}
\widetilde{\kappa} = H^2\left( \sum_{x_j \in \omega_i}  \kappa|\nabla \chi_j|^2\right).
\end{equation*}

Here we have to consider the Rayleigh quotient $$\mathcal{Q}_{ms}(v):=\frac{\int_{\omega_i}  \kappa  |\nabla v|^2}{\int_{\omega_i} \widetilde{\kappa}| v  |^2},$$
with $v\in P^1(\omega_i)$, and define $ I_0^{\omega_i}v $ as projection on low modes. Additional modes ``complement'' the initial space spanned by the partition of unity used so that the resulting coarse space leads to robust  methods with minimal dimension coarse spaces \cite{ge09_1reduceddim}.

If we consider the two-level method with the (multiscale) spectral coarse space presented before, then
\begin{equation*}
\mbox{cond}(M^{-1}A)\leq C(1+(H/\delta)^2),
\end{equation*}
where $C$ is independent of the contrast if enough eigenfunctions in each node neighborhood are selected for the construction of the coarse spaces. The constant  $C$  and the resulting coarse-space dimension depend on the partition of unity (initial coarse-grid representation) used.

\subsection{Abstract eigenvalue problems}\label{secabs}
 We consider an abstract variational problem, where the global bilinear form is obtained by
assembling local bilinear forms. That is
$$a(u,v)=\sum_K a_K(R_Ku,R_Kv),$$ where
$a_K(u,v)$ is a bilinear form acting on functions with supports
being the coarse block $K$. Define the  subdomain bilinear form $a_{\omega_i}(u,v)=\sum_{K\subset \omega_i}a_K(u,v)$.
We consider the abstract problem\\
\centerline{$\displaystyle a(u,v)=F(v) \quad \mbox{ for all } v\in V.$}
We
introduce $\{\chi_j\}$, a partition of unity subordinated to coarse-mesh blocks and
$\{\xi_i\}$ a partition of unity subordinated to overlapping decomposition
(not necessarily related in this subsection). We also define  the
``Mass'' bilinear form (or energy of cut-off) $m_{\omega_i}$ and
the Rayleigh quotient  $\mathcal{Q}_{abs}$ by
 \[\displaystyle m_{\omega_i} (v,v)  :=  \sum_{j\in \omega_i} a(\xi_i\chi_j v,\xi_i\chi_j v)\quad \mbox{ and } \quad \mathcal{Q}_{abs}(v):=\displaystyle \frac{a_{\omega_i} (v,v) }{m_{\omega_i} (v,v) }.
\]

For the Darcy problem,  we have $$ m_{\omega_i} (v,v)
= \sum_{j\in \omega_i} \int_{\omega_i} \kappa |\nabla(\xi_i\chi_j v )|^2
\preceq \int_{\omega_i}\widetilde{ \kappa} |v|^2.$$
The same analysis can be done by replacing the partition of unity functions by  a partition of the degree of freedom (PDoF); see \cite{Efen_GVass_11}. Let  $\{\chi_j\}$ be PDoF subordinated to coarse mesh neighborhood and
$\{\xi_i\}$ be PDoF subordianted to overlapping decomposition. We define the cut-off bilinear form and quotient,
\[\displaystyle m_{\omega_i} (v,v)  :=  \sum_{j\in \omega_i} a(\xi_i\chi_j v,\xi_i\chi_j v) \quad \mbox{ and } \quad \mathcal{Q}_{abs2}(v):=\displaystyle \frac{a_{\omega_i} (v,v) }{m_{\omega_i} (v,v) }.\]
The previous construction allows applying the same design recursively and therefore to use of the same ideas in a multilevel method; see \cite{egv2011,Efen_GVass_11}.

 \subsection{Generalized Multiscale Finite Element Method (GMsFEM) eigenvalue problem}

 We can consider the
 Rayleigh quotients presented before only in a suitable subspace
that allows a good approximation of low modes. We call these subspaces the  snapshot spaces. Denote by $W_i$ the
snapshot space corresponding to subdomain  $\omega_i$, then we consider the Rayleigh quotient,
$$\displaystyle \mathcal{Q}_{gm}(v):= \frac{a_{\omega_i} (v,v) }{m_{\omega_i} (v,v) } \quad
\mbox{with  $v\in W_i.$}$$
The snapshot space can be obtained  by dimension reduction techniques or similar computations; see \cite{egh12, calo2016randomized}. For example, we can consider the following simple example.
In each subdomain $\omega_i$, $i=1,\dots, N_S$:\\

\begin{enumerate}
    \item Generate forcing terms $f_1,f_2,\dots,f_M$ randomly ($\int _{\omega_i}f_\ell=0$);
    \item Compute the local solutions
$ -\mbox{div}(\kappa \nabla u_\ell )=f_\ell $
with homogeneous Neumann boundary condition; 
\item Generate $W_i=\mbox{span}\{ u_\ell \}\cup\{1\}$;
\item  Consider $\mathcal{Q}_{gm}$ with $W_i$ in (3) and compute important modes. 
\end{enumerate}

\section{Numerical results} \label{sec:exp}
\begin{figure}
     \centering
     \begin{subfigure}[b]{0.45\textwidth}
         \centering
         \includegraphics[width=\textwidth]{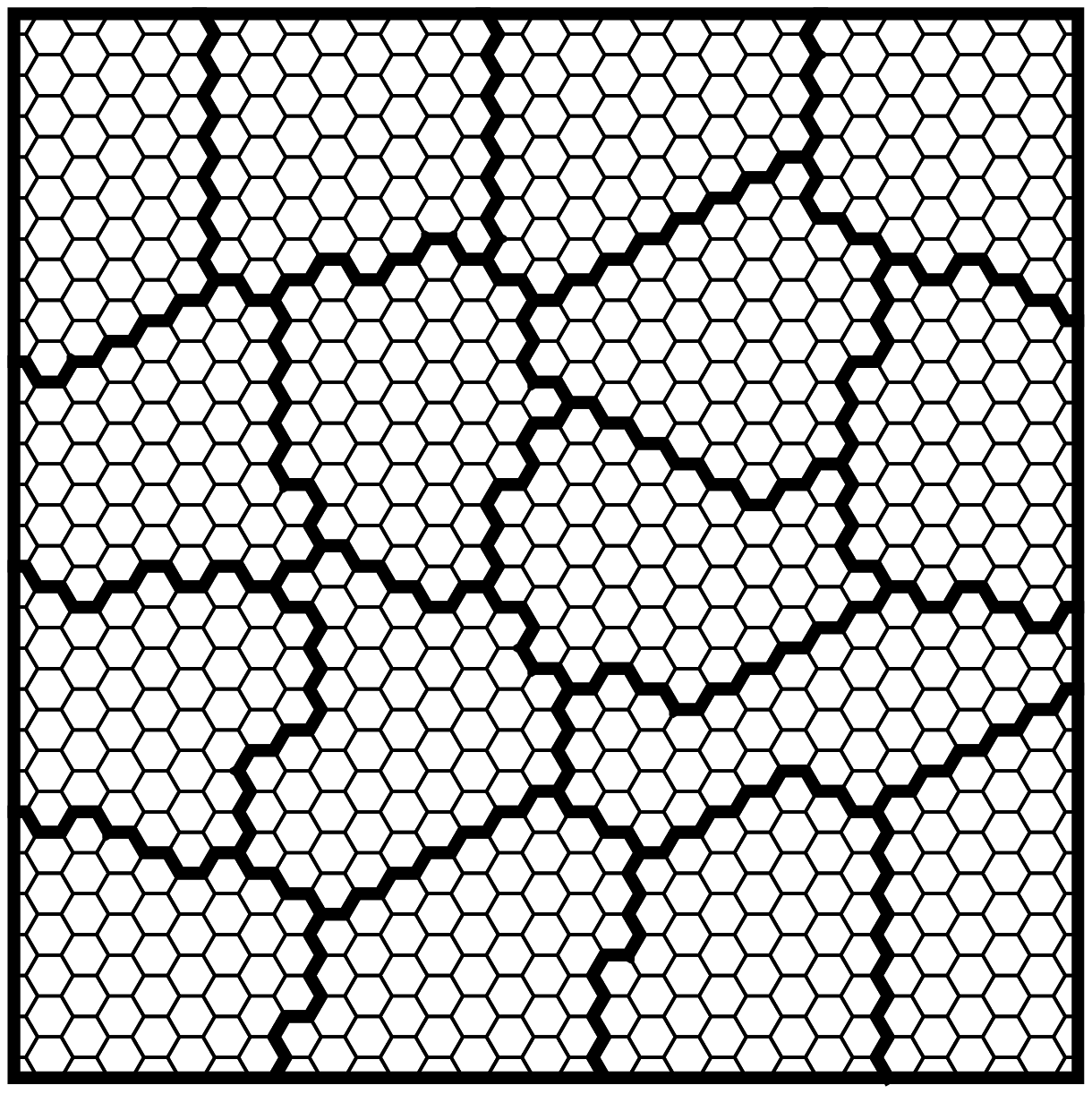}
     \end{subfigure}
     \hfill
     \begin{subfigure}[b]{0.45\textwidth}
         \centering
         \includegraphics[width=\textwidth]{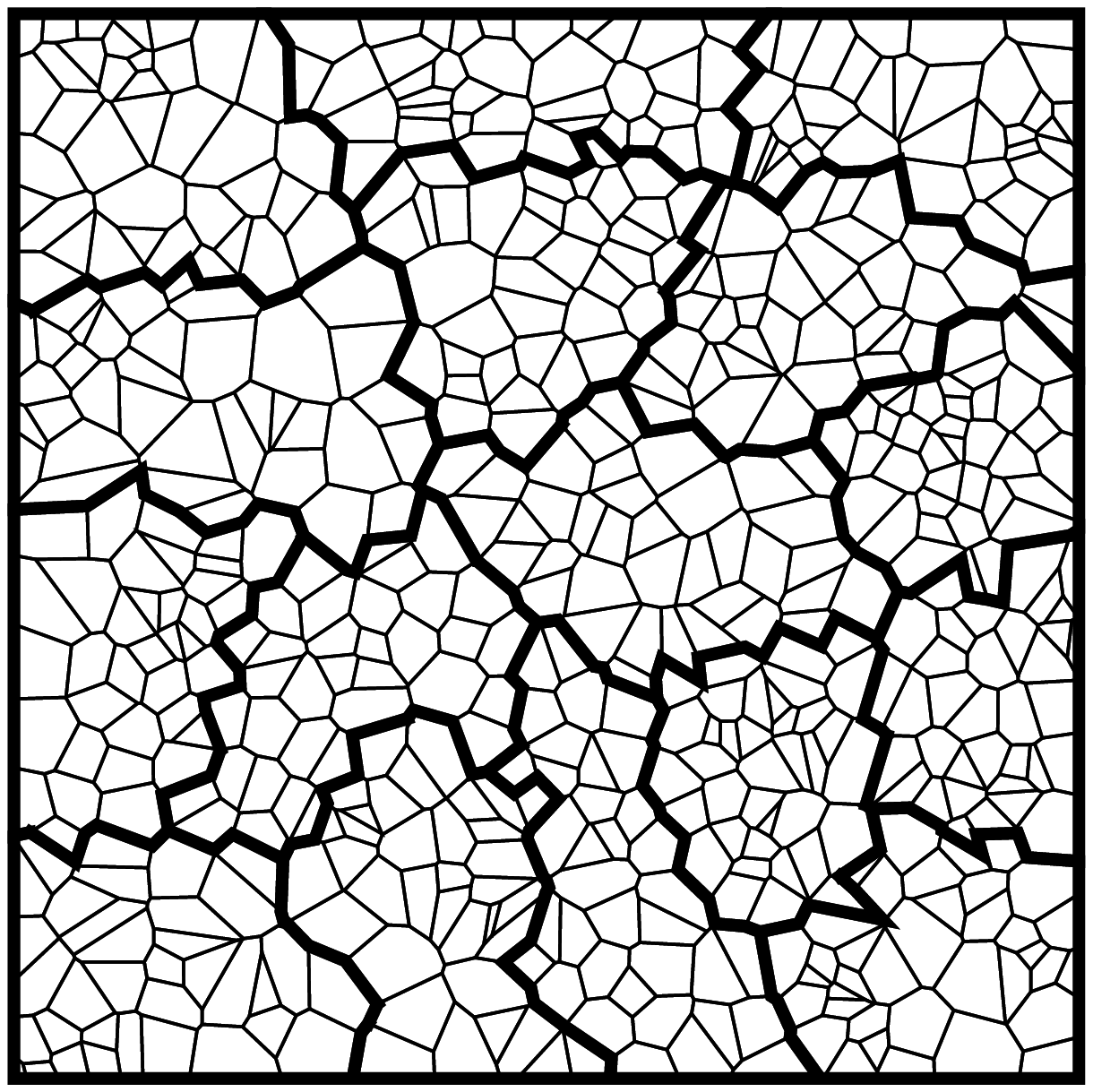}
     \end{subfigure}
        \caption{Fine mesh and 16 subdomains (thick black lines). We present (left) an hexagonal mesh and (right) a Voronoi mesh, with METIS subdomains.}
        \label{fig:meshes}
\end{figure}

In this section we present some numerical results that confirm the behavior of our algorithm for different types of elements, non-overlapping subdomain decompositions (including irregular boundaries), and different high-contrast coefficients.  
We verify that the proposed preconditioner is robust with respect to the high-contrast, multiscale variations of the coefficients and also the presence of irregular boundaries of subdomains. 
We discretize the equation 
$-\mbox{div}(\kappa\nabla u)=0$ with 
homogeneous Dirichlet boundary condition on $\partial D$ using VEM as described in Section \ref{sec:vem}, and then solve the resulting linear system using the preconditioned conjugate gradient method to a relative residual tolerance of $10^{-6}$. We show different meshes and subdomains in Figure \ref{fig:meshes}, and present the behavior of our method as we vary the contrast. We track the condition number, number of iteration until convergence and the dimension of the constructed coarse spaces since these are the main indicators of the robustness of the preconditioner. 

\subsection{Triangular mesh}
\begin{figure}
     \centering
     \begin{subfigure}[b]{0.45\textwidth}
         \centering
         \includegraphics[width=\textwidth]{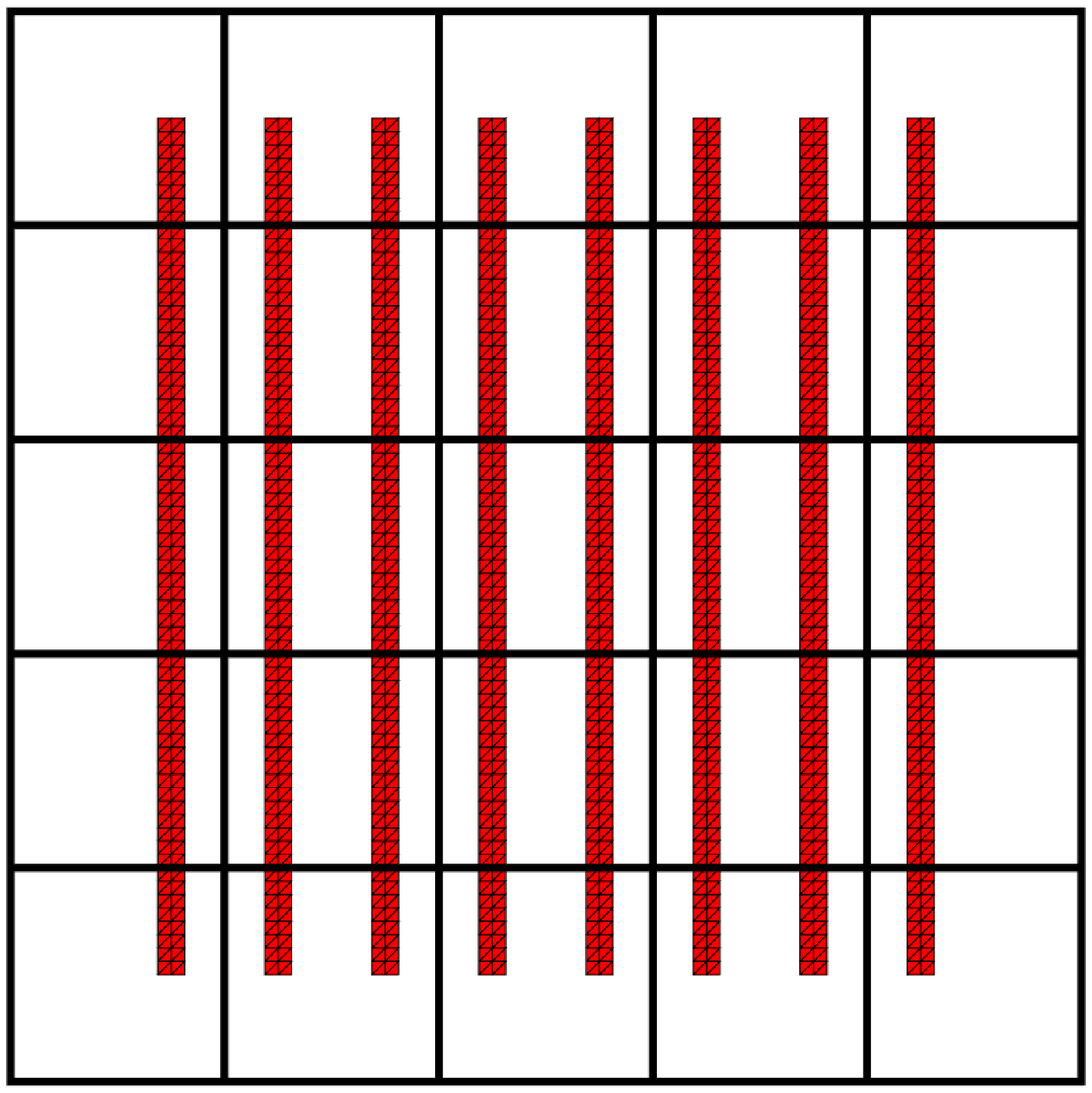}
     \end{subfigure}
     \hfill
     \begin{subfigure}[b]{0.45\textwidth}
         \centering
         \includegraphics[width=\textwidth]{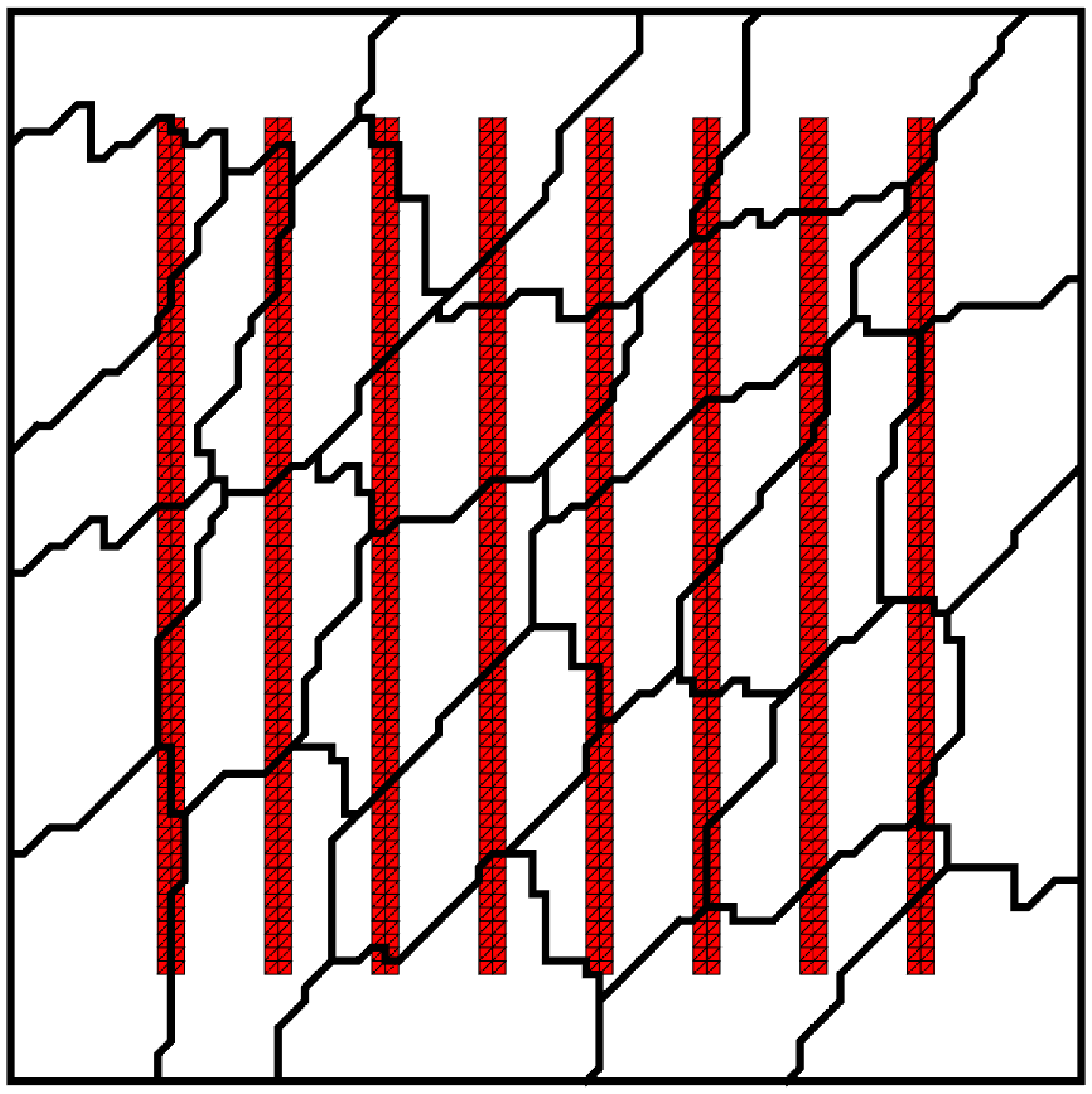}
     \end{subfigure}
        \caption{Function $\kappa$ used in the experiments. We use $\kappa = \eta \in\{1,10^2,10^4,10^6\}$ for red elements, and $\eta=1$ in the background. We show (left) square and (right) METIS subdomains (obtained from a triangular fine mesh); see results in Table \ref{tab:perm3}.}
        \label{fig:tria_channels}
\end{figure}

For comparison with previous studies, we first include results for a triangular mesh with square and METIS subdomains. We consider the methods given in \cite{Widlund_2008} (based on discrete harmonic extensions) and \cite{calvo2018virtual} (based on $L^2-$projections of gradients over polynomial spaces of degree $k$ for virtual element spaces) for the construction of the partition of unity; see results in Table \ref{tab:perm3}. The coefficient $\kappa$ is shown in Figure \ref{fig:tria_channels} for both types of subdomain.

\begin{table}[htb]
\centering
\footnotesize
\begin{tabular}{|c|c||c|c|c||c|c|c|}\hline
\multirow{2}{*}{Subdomains} & \multirow{2}{*}{$\eta$} & \multicolumn{3}{c||}{Non adaptive, harmonic} &\multicolumn{3}{c|}{Adaptive, $k=2$/harmonic}\\ \cline{3-8}
& & Cond &  Iter& dimV0 & Cond &  Iter& dimV0   \\\hline\hline
\multirow{4}{*}{Squares} 
& $1e0$ &17.3  &22 & 36 & 23.8/23.8&21/21 & 4  \\ 
& $1e2$ &44.5  &36 & 36 & 26.8/21.0&28/26 & 56 \\
& $1e4$ &4827  &87 & 36 &  5.0/6.3&18/18 & 96 \\ 
& $1e6$ &1.7e6 &148& 36 &  5.3/5.3&19/18 & 96 \\ 
 \hline\hline
\multirow{4}{*}{METIS}
& $1e0$ &17.8  &29 & 52 & 25.4/25.3&34/34 & 16  \\ 
& $1e2$ &31.3  &45 & 52 & 28.4/12.7&43/27 & 96  \\
& $1e4$ &3741  &167& 52 &  6.9/6.0 &23/22 & 190 \\ 
& $1e6$ &3.0e5 &342& 52 &  7.2/6.0 &25/24 & 190 \\ \hline 
\end{tabular}
\caption{Number of iterations (Iter) until convergence 
of the PCG and condition number (Cond), for different values of the contrast $\eta$, with $\kappa$ as shown in Figure \ref{fig:tria_channels}, for a triangular mesh with $12800$ elements, $25$ subdomains, $H/h\approx 16$, $H/\delta \approx 4$. }
\label{tab:perm3}
\end{table}

\begin{figure}
     \centering
     \begin{subfigure}[b]{0.45\textwidth}
         \centering
         \includegraphics[width=\textwidth]{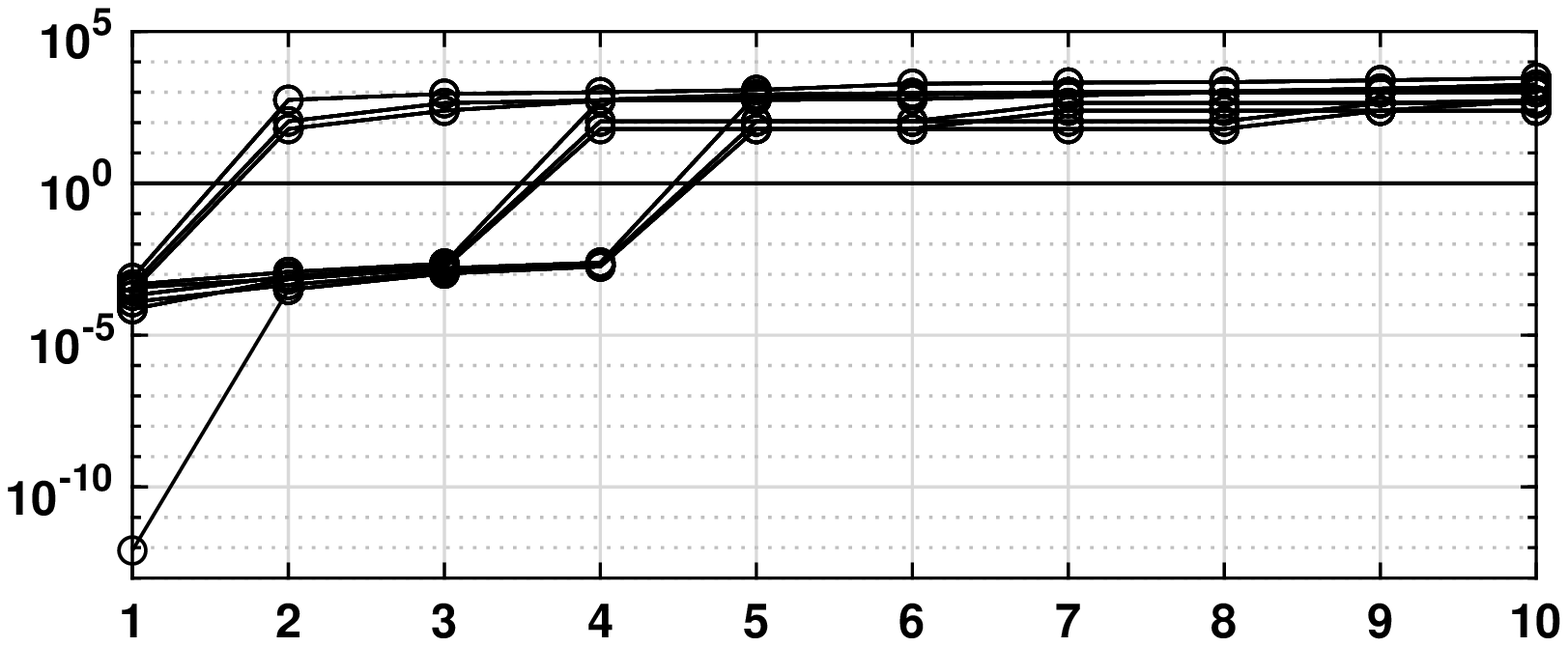}
     \end{subfigure}
     \hfill
     \begin{subfigure}[b]{0.45\textwidth}
         \centering
         \includegraphics[width=\textwidth]{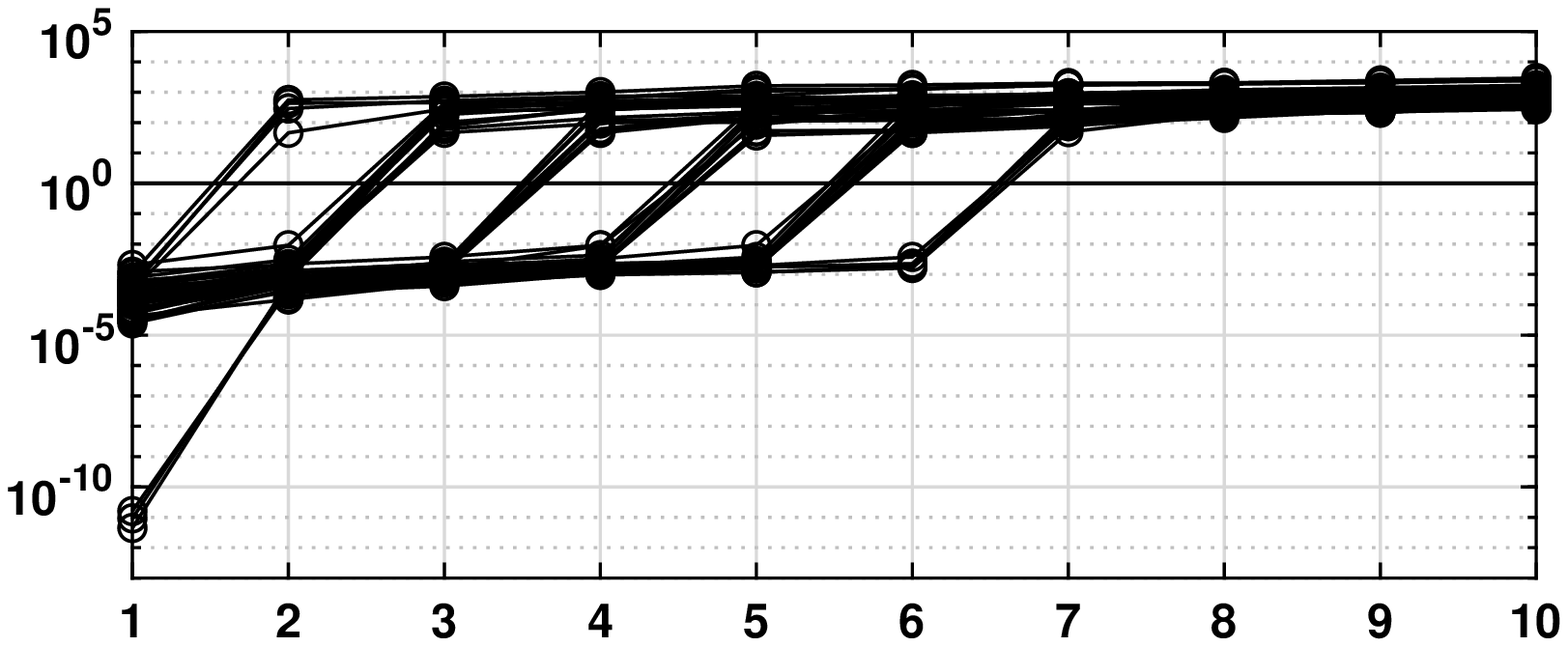}
     \end{subfigure}
        \caption{Smaller eigenvalues (circles) for (left) square and (right) METIS subdomains with a triangular mesh with 12800 elements and $\eta = 10^6$, for which a maximum number of 4 and 6 coarse eigenfuntions per subdomain are required, respectively. We include in the adaptive coarse space all eigenvalues smaller than 1 (black horizontal line); see Table \ref{tab:perm3}.}
        \label{fig:tria_eigs}
\end{figure}

\begin{figure}
         \centering
         \includegraphics[width=\textwidth]{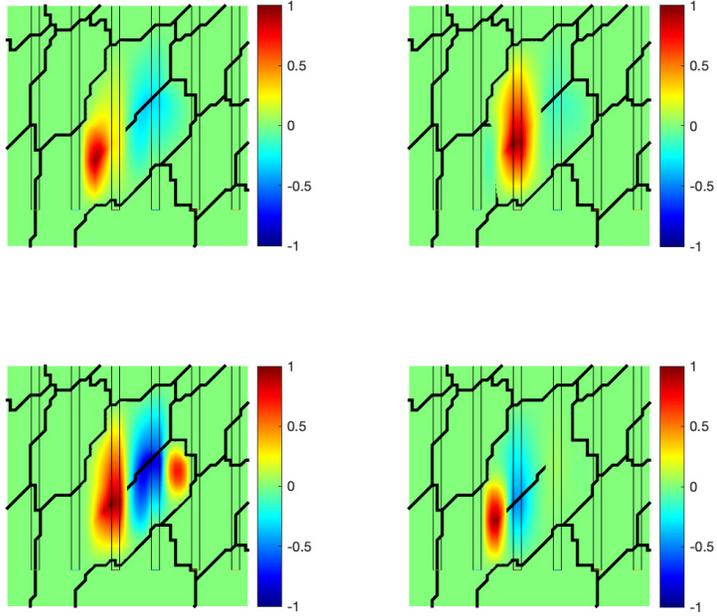}
        \caption{Four eigenfunctions associated to a subdomain vertex for a triangular mesh and METIS subdomains. Thick black lines correspond to the subdomains and thin black lines correspond to the boundary of the channels.}
        \label{fig:eigenfunctions}
\end{figure}

In Figure \ref{fig:tria_eigs} we present the smaller eigenvalues for all the subdomains of the decomposition. We include in the coarse space all eigenfunctions associated to eigenvalues smaller than 1; this implies at most 4 and 6 eigenfunctions per subdomain, for square and METIS subdomains, respectively. Finally, in Figure \ref{fig:eigenfunctions} we include the corresponding eigenfunctions for a subdomain vertex.

\subsection{Hexagonal and Voronoi meshes}

\begin{figure}
     \centering
     \begin{subfigure}[b]{0.45\textwidth}
         \centering
         \includegraphics[width=\textwidth]{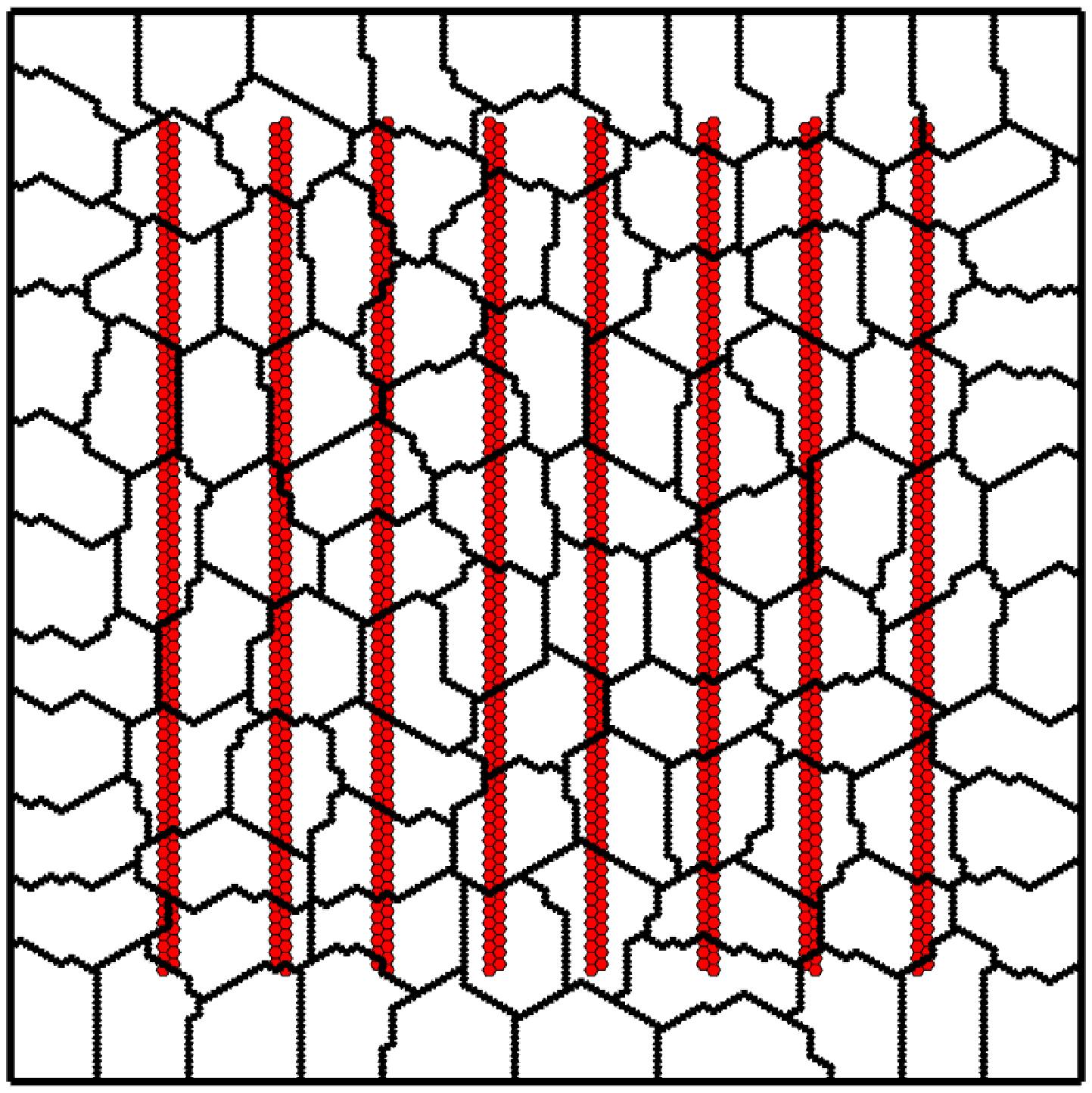}
     \end{subfigure}
     \hfill
     \begin{subfigure}[b]{0.45\textwidth}
         \centering
         \includegraphics[width=\textwidth]{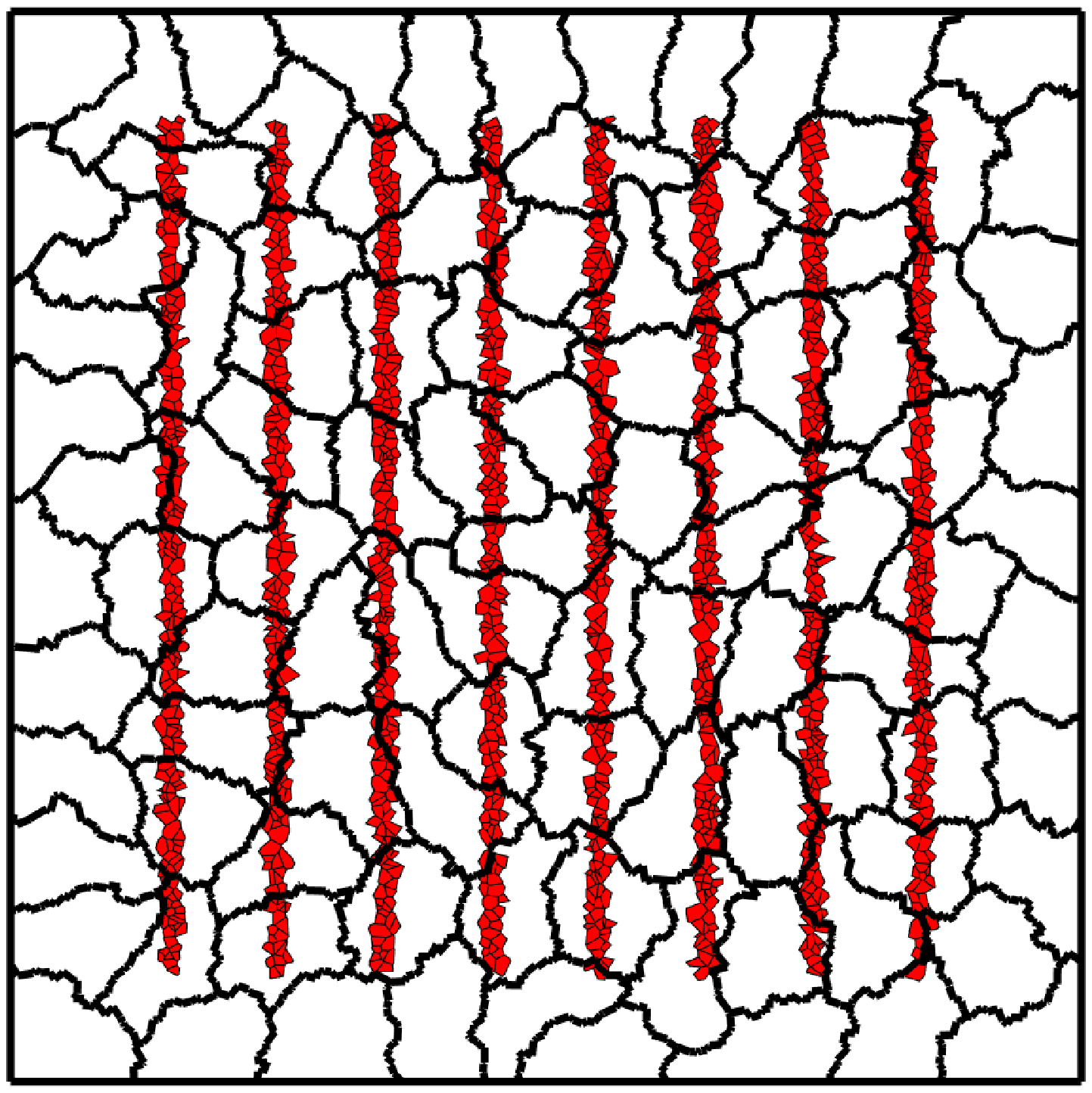}
     \end{subfigure}
        \caption{Function $\kappa$. We use $\kappa = \eta \in\{1,10^2,10^4,10^6\}$ for red elements, and $\eta=1$ in the background. We show (left) hexagonal and (right) Voronoi-type meshes with METIS subdomains.}
        \label{fig:hexa_channels}
\end{figure}

We now consider a hexagonal and a Voronoi mesh, both with METIS subdomains; see Figure \ref{fig:meshes}. We use a similar distribution of high-contrast channels; see Figure \ref{fig:hexa_channels}. We report the number of iterations and condition number estimates in Tables \ref{tab:hexas} and \ref{tab:perm3two}. We observe that the non-adaptive method (as in \cite{Widlund_2008}) deteriorates as $\eta$ increases. By enriching the coarse space as discussed in Section 4.3, we obtain a robust method with respect to the contrast.

\begin{table}[htb]
\centering
\footnotesize
\begin{tabular}{|c||c|c|c||c|c|c||c|c|c|}\hline 
\multirow{2}{*}{$\eta$} & \multicolumn{3}{c||}{Non adaptive, harmonic}  & \multicolumn{3}{c||}{Adaptive, harmonic} & \multicolumn{3}{c|}{Adaptive, $k=2$} \\ \cline{2-10}
 & Cond &  Iter & dimV0 & Cond &  Iter& dimV0 & Cond &  Iter& dimV0  \\\hline
 $1e0$ & 29.6 &  35& 202& 52.1& 48& 96 & 47.7& 48& 96\\
 $1e2$ & 34.3 &  45& 202& 20.8& 37& 191& 52.5& 61& 191\\
 $1e4$ & 1027 & 129& 202& 11.2& 30& 335& 12.8& 34& 335\\
 $1e6$ & 1.0e5& 236& 202& 11.8& 34& 335& 13.6& 37& 335\\
 \hline
\end{tabular}
\caption{Number of iterations (Iter) until convergence 
of the PCG and condition number (Cond), for different values of the contrast $\eta$, with $\kappa$ as shown in Figure \ref{fig:hexa_channels}, for an hexagonal mesh with 9699 elements, 100 subdomains, $H/h\approx 20$, $\delta\approx 2h$.}
\label{tab:hexas}
\end{table}

\begin{table}[htb]
\centering
\footnotesize
\begin{tabular}{|c||c|c|c||c|c|c||c|c|c|}\hline 
\multirow{2}{*}{$\eta$} & \multicolumn{3}{c||}{Non adaptive, harmonic}  & \multicolumn{3}{c||}{Adaptive, harmonic} & \multicolumn{3}{c|}{Adaptive, $k=2$} \\ \cline{2-10}
 & Cond &  Iter & dimV0 & Cond &  Iter& dimV0 & Cond &  Iter& dimV0  \\\hline
 $1e0$ &  43.6& 42& 202& 62.7& 54&  97& 68.3& 55 & 97 \\ 
 $1e2$ &  45.5& 48& 202& 21.1& 37& 199& 63.5& 53 & 199\\
 $1e4$ &  1448&130& 202& 13.3& 34& 382& 17.1& 39 & 382\\ 
 $1e6$ & 1.2e5&246& 202& 13.6& 37& 382& 15.7& 42 & 382\\ 
 \hline
\end{tabular}
\caption{Number of iterations (Iter) until convergence 
of the PCG and condition number (Cond), for different values of the contrast $\eta$, with $\kappa$ as shown in Figure \ref{fig:hexa_channels}, for a Voronoi mesh with 12325 elements, 100 subdomains, $H\approx 34$, $h\approx 0.136$, $\delta\approx 2h$.}
\label{tab:perm3two}
\end{table}

\section{Final comments}\label{sec:conc}
In conclusion, this paper presented an investigation of the efficient solution of PDEs using DDM for problems with high-contrast multiscale coefficients and irregular subdomains. The focus was on addressing the challenges posed by the simultaneous presence of these difficulties, which often arise in real-world applications such as subsurface flow modeling in porous media. The proposed methods required the adaptation of the coarse space and the construction of appropriate partition of unity functions to handle the variations in coefficient and irregular subdomains.
The proposed approach aimed to achieve robustness in DDM performance, minimizing the deterioration caused by high-contrast coefficients, irregular subdomains, and multiscale variations. By combining the constructions of coarse spaces for high-contrast multiscale problems and irregular subdomains, the paper offered a method that addresses these challenges simultaneously.

Overall, the results presented in this paper demonstrated the effectiveness of the proposed approach in efficiently computing solutions for PDEs with complex coefficient structures and irregular subdomains. The numerical experiments provided evidence of the method's robustness and accuracy. The findings contribute to advancing the field of computational mathematics and engineering, providing valuable insights for the solution of practical problems involving the numerical approximation of PDEs in two or three dimensions.\\

\section*{Acknowledgments} 
The first author gratefully acknowledges the
institutional support for the project C1228 subscribed to the Vice-Rectory for Research, University of Costa Rica. The authors express their gratitude to the organizers of the special session titled ``Applied Math and Computational Methods and Analysis across the Americas" at the Mathematical Congress of the Americas 2021, which took place online. This session provided a valuable platform for the authors to convene, engage in fruitful discussions, and exchange ideas related to the subject matter of the paper.

\bibliographystyle{plain}
\bibliography{final}

\end{document}